\documentclass{eptcs-modified}
 \providecommand{\event}{} % Name of the event you are submitting to
 \usepackage{breakurl}

\usepackage{graphicx}
\usepackage{epsfig}
\usepackage{times}
\usepackage{amsthm}
\usepackage{mathptmx}
\usepackage{amssymb,amsfonts,mathrsfs,latexsym,textcomp}
\usepackage{amsmath}
\usepackage{dcolumn}
\usepackage{multicol}
\usepackage{amscd}

\usepackage{epstopdf}

\long\def\greybox#1{%
    \newbox\contentbox%
    \newbox\bkgdbox%
    \setbox\contentbox\hbox to \hsize{%
        \vtop{
            \kern\columnsep
            \hbox to \hsize{%
                \kern\columnsep%
                \advance\hsize by -2\columnsep%
                \setlength{\textwidth}{\hsize}%
                \vbox{
                    \parskip=\baselineskip
                    \parindent=0bp
                    #1
                }%
                \kern\columnsep%
            }%
            \kern\columnsep%
        }%
    }%
    \setbox\bkgdbox\vbox{
        \pdfliteral{0.85 0.85 0.85 rg}
        \hrule width  \wd\contentbox %
               height \ht\contentbox %
               depth  \dp\contentbox
        \pdfliteral{0 0 0 rg}
    }%
    \wd\bkgdbox=0bp%
    \vbox{\hbox to \hsize{\box\bkgdbox\box\contentbox}}%
    \vskip\baselineskip%
}

\title{On the topological aspects of the theory of represented spaces\thanks{Prior versions of this manuscript were titled \emph{Compactness and Separation for Represented Spaces} and \emph{A new introduction to the theory of represented spaces}.}}
\author{
Arno Pauly
\institute{Clare College\\ University of Cambridge, United Kingdom}
\email{Arno.Pauly@cl.cam.ac.uk}
}
\def\titlerunning{Theory of Represented Spaces}
\def\authorrunning{A. Pauly}

\begin{document}
 \theoremstyle{definition}
 \newtheorem{theorem}{Theorem}
\newtheorem{definition}[theorem]{Definition}
\newtheorem{problem}[theorem]{Problem}
\newtheorem{assumption}[theorem]{Assumption}
\newtheorem{corollary}[theorem]{Corollary}
\newtheorem{proposition}[theorem]{Proposition}
\newtheorem{lemma}[theorem]{Lemma}
\newtheorem{observation}[theorem]{Observation}
\newtheorem{fact}[theorem]{Fact}
\newtheorem{question}[theorem]{Open Question}
\newtheorem{conjecture}[theorem]{Conjecture}
\newcommand{\dom}{\operatorname{dom}}
\newcommand{\id}{\textnormal{id}}
\newcommand{\Cantor}{{\{0, 1\}^\mathbb{N}}}
\newcommand{\Baire}{{\mathbb{N}^\mathbb{N}}}
\newcommand{\Lev}{\textnormal{Lev}}
\newcommand{\hide}[1]{}
\newcommand{\mto}{\rightrightarrows}
\newcommand{\uint}{{[0, 1]}}
\newcommand{\bft}{\mathrm{BFT}}
\newcommand{\lbft}{\textnormal{Linear-}\mathrm{BFT}}
\newcommand{\pbft}{\textnormal{Poly-}\mathrm{BFT}}
\newcommand{\sbft}{\textnormal{Smooth-}\mathrm{BFT}}
\newcommand{\ivt}{\mathrm{IVT}}
\newcommand{\cc}{\textrm{CC}}
\newcommand{\lpo}{\textrm{LPO}}
\newcommand{\llpo}{\textrm{LLPO}}
\newcommand{\aou}{AoU}
\newcommand{\Ctwo}{C_{\{0, 1\}}}
\newcommand{\name}[1]{\textsc{#1}}
\newcommand{\C}{\textrm{C}}
\newcommand{\ic}[1]{\textrm{C}_{\sharp #1}}
\newcommand{\xc}[1]{\textrm{XC}_{#1}}
\newcommand{\me}{\name{P}.~}
\newcommand{\etal}{et al.~}
\newcommand{\eval}{\operatorname{eval}}
\newcommand{\Sierp}{Sierpi\'nski }
\newcommand{\isempty}{\operatorname{IsEmpty}}
\newcommand{\nonempty}{\operatorname{IsNonEmpty}}
\newcommand{\sat}[1]{\mathalpha{\uparrow}#1}

\maketitle

\begin{abstract}
Represented spaces form the general setting for the study of computability derived from Turing machines. As such, they are the basic entities for endeavors such as computable analysis or computable measure theory. The theory of represented spaces is well-known to exhibit a strong topological flavour. We present an abstract and very succinct introduction to the field; drawing heavily on prior work by \name{Escard\'o}, \name{Schr\"oder}, and others.

Central aspects of the theory are function spaces and various spaces of subsets derived from other represented spaces, and -- closely linked to these -- properties of represented spaces such as compactness, overtness and separation principles. Both the derived spaces and the properties are introduced by demanding the computability of certain mappings, and it is demonstrated that typically various interesting mappings induce the same property.
\end{abstract}

\section{Introduction}
Just as numberings provide a tool to transfer computability from $\mathbb{N}$ (or $\{0,1\}^*$) to all sorts of countable structures; representations provide a means to introduce computability on structures of the cardinality of the continuum based on computability on Cantor space $\Cantor$. This is of course essential for computable analysis \cite{weihrauchd}, dealing with spaces such as $\mathbb{R}$, $\mathcal{C}(\mathbb{R}^n, \mathbb{R}^m)$, $\mathcal{C}^k(\mathbb{R}^n, \mathbb{R}^m)$ or general Hilbert spaces \cite{brattka9}. Computable measure theory (e.g.~\cite{weihrauche,schroder2,paulymeasurement,hoyrup,collins4}), or computability on the set of countable ordinals \cite{zhenhao,pauly-ordinals} likewise rely on representations as foundation.

Essentially, by equipping a set with a representation, we arrive at a \emph{represented space} -- to some extent\footnote{The notion of a multi-representation \cite{schroder3,weihrauchg} goes beyond representations. A multi-representation of a set $X$ is a right-total relation $\delta \subseteq \Cantor \times X$, relating codes to the encoded elements. The distinction to representation is that a code can stand for more than one element here.

Several common examples of multi-representations have the additional feature that $\exists p \in \Cantor \delta(p, x) \wedge \delta(p, y)$ defines an equivalence relation on $X$ -- hence, we can conceive of the multi-representation as an ordinary representation of the induced equivalence classes. Often we even have a canonization operation available; such as the saturation used in Section \ref{sec:compactness} and the closure used in Section \ref{sec:overt} to obtain represented spaces.

A very different example (communicated to the author by \name{Bauer}) is the multirepresentation $\mathalpha{\in} : \Cantor \mto \left (\mathcal{P}(\Cantor) \setminus \{\emptyset\}\right)$. Examples of this kind, however, are beyond the scope of the paper.} the most general structure carrying a notion of computability that is derived from Turing machines\footnote{Algebraic computation models \cite{blum2} as abstract generalizations of register machines are a completely different thing.}. Any represented space provides in a canonic way also function spaces and spaces of certain subsets. While the usefulness of these derived spaces will vary depending on the area of application, their development can be done while abstracting away from the specifics of the original represented space -- in fact, such details have traditionally obfuscated the core proof ideas. Closely linked to the various subset spaces, also a few properties of spaces defined in terms of computability of certain maps are studied. Among these, special attention shall be drawn to admissibility, which has often been advocated as a criterion for well-behaved representations in computable analysis \cite{weihrauchd, schroder}.

\section{Connections to the literature}
It is well-known that many representations in computable analysis can be characterized by extremal properties -- the prototypic example being that the standard representation of the space of continuous functions carries exactly as much information as needed to make function evaluation computable. In fact, this condition not only characterizes the representation, but also the set itself: Whenever a function is contained in a function space admitting computable function evaluation, it is computable w.r.t.~some oracle -- its name -- hence, it is continuous. Consequently, rather than speaking merely about representations (of fixed sets) being characterized such, we should consider \emph{represented spaces} -- the combination of a set and a representation of it -- as the fundamental objects of the characterization results\footnote{Hence the chance of the term of choice from \name{Kreitz}' and \name{Weihrauch}'s \emph{Theory of Representations} \cite{kreitz} to the present paper's \emph{Theory of Represented Spaces}.}.

In many characterization results in the literature the representation is defined explicitly first, then the extremal property is proven. However, as a consequence of the UTM theorem, this is unnecessary: The condition itself defines the represented space up to computable isomorphisms. In fact, specifying details of standard representations beyond their characterizing property often complicates proofs of basic results, and can obfuscate the algorithmic ideas.

Developing the basic theory of represented spaces based on extremal characterizations amounts to an instantiation of \name{Escard\'o}'s \emph{synthetic topology} \cite{escardo, escardo2} with the category of represented spaces. As this category is particularly well-behaved, we obtain a very concise picture, and can prove a few more results on compactness and separation than available in the generic setting. Similarly, there are many parallels to \name{Taylor}'s \emph{abstract stone duality} \cite{taylor, taylor2}, in comparison we again sacrifice generality for simplicity of presentation and strength of some results.

Most definitions and many results about them were already formulated by \name{Schr\"oder} \cite{schroder5} before the advent of synthetic topology though. In this development, spaces were always required to satisfy some form of admissibility\footnote{The form of admissibility presented in this paper corresponds to admissibility w.r.t.~some topology. Schr\"oder also studied admissibility w.r.t~some variants of limit spaces.}. For our purposes, these restrictions are mostly immaterial, and will be dropped. Note that admissibility (w.r.t.~a topology) turns out to be a feature of synthetic topology in general, demarking those spaces actually fully understandable in terms of their (internal) topology. For represented spaces, this can be seen as (slightly) generalizing \cite{schroder}. The PhD thesis of \name{Lietz} \cite{lietz} also contains relevant results on admissibility, in particular \cite[Theorem 3.2.7.]{lietz}. A direct application of the framework of synthethic topology to represented spaces was performed by \name{Collins} \cite{collins3,collins5} in 2010, several of the result presented here are already present in \cite{collins5}.

Effective compactness has been studied in \cite{weihrauchj} by \name{Brattka} and \name{Weihrauch}, and in \cite{presser} by \name{Brattka} and \name{Presser}. \cite{grubba} by \name{Grubba} and \name{Weihrauch} considers representations of closed, open and compact sets characterized by extremal properties for the restricted case of countably based spaces. Effective topological separation was considered in \cite{weihrauchh,weihrauchm} by \name{Weihrauch}.

A warning is due regarding a discrepancy in nomenclature: The computable elements of the space of closed subsets introduced here (i.e. the computable closed sets) are called co-c.e. closed sets by \name{Weihrauch} and others. Their c.e. closed sets are referred to as \emph{overt} here (a term coined by \name{Taylor}, see Section \ref{sec:overt}), as they as a space lack the structure expected from closed sets. Subsequently, the sets called computabl{\bf y} closed in \name{Weihrauch}'s terminology are called computable closed and overt here.

The observation that the c.e. closed sets lack the closure properties of closed sets spawned various investigations into more restricted setting remedying the issue. An example of this is \name{Ziegler}'s work on representations for regular closed subsets of a Euclidean space \cite{ziegler8}. A second example, somewhat further removed from representations, is the observation that computable closed sets (i.e.~co-c.e. closed) satisfying some topological criterion have to be computably overt (i.e.~c.e. closed) as made by \name{Miller} and \name{Iljazovi\'c} \cite{miller,miller3,iljazovic,iljazovic2}.

As the attribution of the definitions and results is non-trivial due to parallel independent developments using slightly different formal frameworks, the discussion of the intellectual pedigree follows separately each section. In addition to work cited there, several observations presumably are folklore.

This paper is originally based on \cite[Chapter 3]{paulyphd}, the PhD thesis of the author.

\section{The foundations}
The central notion is the \emph{represented space}, which is a pair $\mathbf{X} = (X, \delta_X)$ of a set $X$ and a partial surjection $\delta_X : \subseteq \Cantor \to X$ (\footnote{Our choice of Cantor space as domain for our representations is inconsequential, often Baire space $\Baire$ is used instead. A recent observation by \name{Cook} and \name{Kawamura} is that using the set of regular functions on natural numbers as the foundation may be useful for complexity theory \cite{kawamura}, again, this has no impact at all on the computability theory presented here.}). A multi-valued function between represented spaces is a multi-valued function between the underlying sets. For $f : \subseteq \mathbf{X} \mto \mathbf{Y}$ and $F : \subseteq \Cantor \to \Cantor$, we call $F$ a realizer of $f$ (notation $F \vdash f$), iff $\delta_Y(F(p)) \in f(\delta_X(p))$ for all $p \in \dom(f\delta_X)$. A map between represented spaces is called computable (continuous), iff it has a computable (continuous) realizer. Similarly, we call a point $x \in \mathbf{X}$ computable, iff there is some computable $p \in \Baire$ with $\delta_\mathbf{X}(p) = x$. We write $\mathbf{X} \cong \mathbf{Y}$ if the two spaces are computably isomorphic, i.e.~if there is a bijection $f : \mathbf{X} \to \mathbf{Y}$ that is computable and has a computable inverse. For our purposes, there is no need to distinguish computably isomorphic spaces.

\greybox{{\bf Warning:} A priori, the notion of a continuous map between represented spaces and a continuous map between topological spaces are distinct and should not be confused!
}
We consider two categories of represented spaces, one equipped with the computable maps, and one equipped with the continuous maps. We call the resulting structure a \emph{category extension} (cf.~\cite{paulysearchproblems,paulysearchproblemscie}), as the former is a subcategory of the latter, and shares its structure (products, coproducts, exponentials) as we shall see next. In general, all our results have two instances, one for the computable and one for the continuous maps, with the proofs being identical. Essentially, the continuous case can be considered as the relativization of the computable one, as a function on Cantor space is continuous iff it is computable relative to some oracle.

\begin{proposition}
\label{prop:categorystructure}
For any two represented spaces $\mathbf{X}$, $\mathbf{Y}$ there are represented spaces $\mathbf{X} \times \mathbf{Y}$, $\mathbf{X} + \mathbf{Y}$ and computable functions $\pi_1 : \mathbf{X} \times \mathbf{Y} \to \mathbf{X}$, $\pi_2 : \mathbf{X} \times \mathbf{Y} \to \mathbf{X}$, $\iota_1 : \mathbf{X} \to \mathbf{X} + \mathbf{Y}$, $\iota_2 : \mathbf{Y} \to \mathbf{X} + \mathbf{Y}$, such that: \begin{enumerate}
\item For any pair of continuous maps $f_1 : \mathbf{Z} \to \mathbf{X}$, $f_2 : \mathbf{Z} \to \mathbf{Y}$ there is a unique continuous map $\langle f_1, f_2\rangle : \mathbf{Z} \to \mathbf{X} \times \mathbf{Y}$ with $f_i = \pi_i \circ \langle f_1,f_2\rangle$.
\item If $f_1$ and $f_2$ are computable, so is $\langle f_1,f_2\rangle$.
\item For any pair of continuous maps $f_1 : \mathbf{X} \to \mathbf{Z}$, $f_2 : \mathbf{Y} \to \mathbf{Z}$ there is a unique continuous map $(f_1 + f_2) : \mathbf{X} + \mathbf{Y} \to \mathbf{Z}$ with $f_i = (f_1 + f_2) \circ \iota_i$.
\item If $f_1$ and $f_2$ are computable, so is $f_1 + f_2$.
\end{enumerate}
The space $\mathbf{X} + \mathbf{Y}$ can be obtained via $(X, \delta_X) + (Y, \delta_Y) = (X \biguplus Y, \delta_{X+Y})$ where $\delta_{X+Y}(0p) = \delta_X(p)$ and $\delta_{X+Y}(1p) = \delta_Y(p)$. The space $\mathbf{X} \times\mathbf{Y}$ can be obtained via $(X, \delta_X) \times (Y, \delta_Y) = (X \times Y, \delta_{X\times Y})$ where $\delta_{X\times Y}(\langle p, q\rangle) = (\delta_X(p), \delta_Y(q))$.
\end{proposition}

\begin{definition}
Given a pair of represented spaces $\mathbf{X} = (X, \delta_X)$ and $\mathbf{Y} = (Y, \delta_Y)$ we define $\mathbf{X} \wedge \mathbf{Y} := (X \cap Y, \delta_X \wedge \delta_Y)$ where $(\delta_X \wedge \delta_Y)(\langle p, q\rangle) = x$ iff $\delta_X(p) = x \wedge \delta_Y(q) = x$.
\end{definition}
Given two represented spaces $\mathbf{X}$, $\mathbf{Y}$ we obtain a third represented space $\mathcal{C}(\mathbf{X}, \mathbf{Y})$ of functions from $X$ to $Y$ by letting $0^n1p$ be a $[\delta_X \to \delta_Y]$-name for $f$, if the $n$-th Turing machine equipped with the oracle $p$ computes a realizer for $f$. As a consequence of the UTM theorem (in the form proven by \name{Weihrauch} \cite{weihrauchk}), $\mathcal{C}(-, -)$ is the exponential in the category of continuous maps between represented spaces, and the evaluation map is even computable.

\begin{proposition}
\label{compana:prop:functionspacesbasics}
Let $\mathbf{X} = (X, \delta_X)$, $\mathbf{Y} = (Y, \delta_Y)$, $\mathbf{Z} = (Z, \delta_Z)$, $\mathbf{U} = (U, \delta_U)$ be represented spaces. Then the following functions are computable:
\begin{enumerate}
\item $\eval : \mathcal{C}(\mathbf{X}, \mathbf{Y}) \times \mathbf{X} \to \mathbf{Y}$ defined by $\eval(f, x) = f(x)$.
\item $\operatorname{curry} : \mathcal{C}(\mathbf{X} \times \mathbf{Y}, \mathbf{Z}) \to \mathcal{C}(\mathbf{X}, \mathcal{C}(\mathbf{Y}, \mathbf{Z}))$ defined by $\operatorname{curry}(f) = x \mapsto (y \mapsto f(x, y))$.
\item $\operatorname{uncurry} : \mathcal{C}(\mathbf{X}, \mathcal{C}(\mathbf{Y}, \mathbf{Z})) \to \mathcal{C}(\mathbf{X} \times \mathbf{Y}, \mathbf{Z})$ defined by $\operatorname{uncurry}(f) = (x, y) \mapsto f(x)(y)$.
\item $\circ : \mathcal{C}(\mathbf{Y}, \mathbf{Z}) \times \mathcal{C}(\mathbf{X}, \mathbf{Y}) \to \mathcal{C}(\mathbf{X}, \mathbf{Z})$, the composition of functions
\item $\times : \mathcal{C}(\mathbf{X}, \mathbf{Y}) \times \mathcal{C}(\mathbf{U}, \mathbf{Z}) \to \mathcal{C}(\mathbf{X} \times \mathbf{U}, \mathbf{Y} \times \mathbf{Z})$
\item $\operatorname{const} : \mathbf{Y} \to \mathcal{C}(\mathbf{X}, \mathbf{Y})$ defined by $\operatorname{const}(y) = (x \mapsto y)$.
\end{enumerate}
\begin{proof}
All items follow from standard arguments on Turing machines; one merely has to verify that the Type-2 semantics are unproblematic.
\begin{enumerate}
\item By definition of $[\delta_X \to \delta_Y]$, if we apply the Turing machine with oracle specified in the $[\delta_X \to \delta_Y]$-name of $f$ to a $\delta_X$-name of $x$, we obtain a $\delta_Y$-name of $f(x)$.
\item From a Turing machine $M$ we can compute a Turing machine $M'$, such that $M'$ on input $p$ and oracle $\langle q, o\rangle$ simulates $M$ on input $\langle p, q\rangle$ and oracle $o$.
\item and vice versa.
\item Composition of Turing machines is computable, and appropriate access to the two oracles can be ensured. As composition of realizers yields a realizer of the composition, this suffices.
\item The execution of two Turing machines in parallel can be simulated by a single one, access to the oracles can be done accordingly. Products of realizers are realizers of products.
\item This is done via a Turing machine which ignores its input and copies the oracle tape to the output tape.
\end{enumerate}
\end{proof}
\end{proposition}

\begin{corollary}
\label{compana:corr:partial}
Let $\mathbf{X}$, $\mathbf{Y}$, $\mathbf{Z}$ be represented spaces. For any $x \in X$, the map $\operatorname{partial}_x : \mathcal{C}(\mathbf{X} \times \mathbf{Y}, \mathbf{Z}) \to \mathcal{C}(\mathbf{Y}, \mathbf{Z})$ defined by $\operatorname{partial}_x(f) = (y \mapsto f(x, y))$ is continuous. If $x$ is computable, then so is $\operatorname{partial}_x$.
\end{corollary}

\subsubsection*{Remarks}
The earliest use of \emph{represented spaces} as entities in their own right seems to be due to \name{Brattka} in 1996 \cite{brattka13}, although the treatment there does not go beyond the introduction of computability before the setting is restricted to (countably based) topological spaces. The notion of a \emph{representation} (with full generality\footnote{In a more restrictive sense, \name{Weihrauch} and \name{sch\"afer} used the term \emph{representation} before \cite{weihrauchl}.}) precedes this by far, going back to \name{Weihrauch} \cite{weihrauchk} and \name{Kreitz} and \name{Weihrauch} \cite{kreitz} in 1985. The category of representation (of represented spaces) was explicitly mentioned by \name{Bauer} in 2002 \cite{bauer3}.

That the study of computability and the study of continuity go hand in hand when it comes to the theory of representations (of represented spaces) has been noted since its beginnings in \cite{kreitz}, and takes a very visible role e.g.~in \name{Weihrauch}'s books \cite{weihrauch, weihrauchd}. That this yields two categories with the same objects and the same structure was observed by \name{Bauer} (in a slightly different setting) in 1998 \cite{bauer5}. The observation that continuity is relativized computability, which is fundamental for this duality, seems to be a folk result\footnote{The difficulty in attributing the result to any specific person is also mentioned in \cite{kapron}.} appearing first in the 1950's.

The operations $+, \times, \wedge$ on representations are all defined already in \cite{kreitz}, their properties as given in Proposition \ref{prop:categorystructure} appear e.g.~in \cite{weihrauchd}. The exponential is introduced in \cite{weihrauchk} (although not named as such). The category of represented spaces is identified as cartesian closed in \cite{bauer3}.

\section{Open and closed sets}
In the following, we will want to make use of two special represented spaces, $\mathbb{N} = (\mathbb{N}, \delta_\mathbb{N})$ and $\mathbb{S} = (\{\bot, \top\}, \delta_\mathbb{S})$. The representation are given by $\delta_\mathbb{N}(0^n10^\mathbb{N}) = n$, $\delta_\mathbb{S}(0^\mathbb{N}) = \bot$ and $\delta_\mathbb{S}(p) = \top$ for $p \neq 0^\mathbb{N}$. It is straightforward to verify that the computability notion for the represented space $\mathbb{N}$ coincides with classical computability over the natural numbers.

The \Sierp space $\mathbb{S}$ in turn allows us to formalize semi-decidability. The computable functions $f : \mathbb{N} \to \mathbb{S}$ are exactly those where $f^{-1}(\{\top\})$ is recursively enumerable (and thus $f^{-1}(\{\bot\})$ co-recursively enumerable). In general, for any represented space $\mathbf{X}$ we obtain two spaces of subsets of $\mathbf{X}$; the space of open sets $\mathcal{O}(\mathbf{X})$ by identifying $f \in \mathcal{C}(\mathbf{X}, \mathbb{S})$ with $f^{-1}(\{\top\})$, and the space of closed sets $\mathcal{A}(\mathbf{X})$ by identifying $f \in \mathcal{C}(\mathbf{X}, \mathbb{S})$ with $f^{-1}(\{\bot\})$. The properties of the spaces of open and closed sets follow from a few particular computable functions on \Sierp space $\mathbb{S}$ and the function space properties in Proposition \ref{compana:prop:functionspacesbasics}.

\begin{proposition}
\label{prop:sierpbasics}
The functions $\wedge, \vee : \mathbb{S} \times \mathbb{S} \to \mathbb{S}$ and $\bigvee : \mathcal{C}(\mathbb{N}, \mathbb{S}) \to \mathbb{S}$ are computable.
\end{proposition}

Here we have $\wedge, \vee : \mathbb{S} \times \mathbb{S} \to \mathbb{S}$ defined by $\wedge(\top, \top) = \top$, $\wedge(x, y) = \bot$ for $(x, y) \neq (\top, \top)$, $\vee(\bot, \bot) = \bot$, $\vee(x, y) = \top$ for $(x, y) \neq (\bot, \bot)$. Moreover, $\bigvee$ is defined by $\bigvee((x_n)_{n \in \mathbb{N}}) = \top$ iff $\exists n_\top \in \mathbb{N}$ s.t.~$x_{n_\top} = \top$, and $\bigvee((x_n)_{n \in \mathbb{N}}) = \bot$ otherwise.

\begin{proof}
A machine realizing $\wedge, \vee$ simply writes $0$s until it reads a non-zero value from one (for $\vee)$ or both (for $\wedge$) input components, then it continues with $1$s. In order to solve $\bigvee$, all $(x_n)$ are investigated simultaneously, while $0$s are written. If a $1$ ever occurs anywhere in the input, the program starts writing $1$s.
\end{proof}

We point out that neither $\neg : \mathbb{S} \to \mathbb{S}$ or $\bigwedge : \mathcal{C}(\mathbb{N}, \mathbb{S}) \to \mathbb{S}$ are computable.

\begin{proposition}
\label{compana:prop:basicsetoperations}
Let $\mathbf{X}$, $\mathbf{Y}$ be represented spaces. Then the following functions are well-defined and computable:
\begin{enumerate}
\item $^C : \mathcal{O}(\mathbf{X}) \to \mathcal{A}(\mathbf{X})$, $^C : \mathcal{A}(\mathbf{X}) \to \mathcal{O}(\mathbf{X})$ mapping a set to its complement
\item $\cup : \mathcal{O}(\mathbf{X}) \times \mathcal{O}(\mathbf{X}) \to \mathcal{O}(\mathbf{X})$, $\cup : \mathcal{A}(\mathbf{X}) \times \mathcal{A}(\mathbf{X}) \to \mathcal{A}(\mathbf{X})$
\item $\cap : \mathcal{O}(\mathbf{X}) \times \mathcal{O}(\mathbf{X}) \to \mathcal{O}(\mathbf{X})$, $\cap : \mathcal{A}(\mathbf{X}) \times \mathcal{A}(\mathbf{X}) \to \mathcal{A}(\mathbf{X})$
\item $\bigcup : \mathcal{C}(\mathbb{N}, \mathcal{O}(\mathbf{X})) \to \mathcal{O}(\mathbf{X})$ mapping a sequence $(U_n)_{n \in \mathbb{N}}$ of open sets to their union $\bigcup_{n \in \mathbb{N}} U_n$
\item $\bigcap : \mathcal{C}(\mathbb{N}, \mathcal{A}(\mathbf{X})) \to \mathcal{A}(\mathbf{X})$ mapping a sequence $(A_n)_{n \in \mathbb{N}}$ of closed sets to their intersection $\bigcap_{n \in \mathbb{N}} A_n$
\item $^{-1} : \mathcal{C}(\mathbf{X}, \mathbf{Y}) \to \mathcal{C}(\mathcal{O}(\mathbf{Y}), \mathcal{O}(\mathbf{X}))$ mapping $f$ to $f^{-1}$ as a set-valued function for open sets
\item $\mathalpha{\in} : \mathbf{X} \times \mathcal{O}(\mathbf{X}) \to \mathbb{S}$ defined by $\mathalpha{\in}(x, U) = \top$, if $x \in U$.
\item $\times : \mathcal{O}(\mathbf{X}) \times \mathcal{O}(\mathbf{Y}) \to \mathcal{O}(\mathbf{X} \times \mathbf{Y})$, $\times : \mathcal{A}(\mathbf{X}) \times \mathcal{A}(\mathbf{Y}) \to \mathcal{A}(\mathbf{X} \times \mathbf{Y})$
\item $\operatorname{Cut} : \mathbf{Y} \times \mathcal{O}(\mathbf{X} \times \mathbf{Y}) \to \mathcal{O}(\mathbf{X})$ mapping $(y, U)$ to $\{x \mid (x, y) \in U\}$
\item $\Pi : \mathcal{C}(\mathbb{N}, \mathcal{A}(\mathbf{X})) \to \mathcal{A}(\mathcal{C}(\mathbb{N}, \mathbf{X}))$, where $\Pi((A_n)_{n \in \mathbb{N}}) = \Pi_{n \in \mathbb{N}} A_n$.
\end{enumerate}
\begin{proof}
\begin{enumerate}
\item By definition both functions are realized by $\id_\Cantor$.
\item Taking into consideration that $\mathcal{O}(\mathbf{X})$, $\mathcal{A}(\mathbf{X})$ may be considered as subspaces of the function space $\mathcal{C}(\mathbf{X}, \mathbb{S})$, we may realize $\cup : \mathcal{O}(\mathbf{X}) \times \mathcal{O}(\mathbf{X}) \to \mathcal{O}(\mathbf{X})$ by composition (Proposition \ref{compana:prop:functionspacesbasics} (4, 6), together with the diagonal) of $\wedge$ and $\times$ from Proposition \ref{compana:prop:functionspacesbasics} (5), likewise the composition of $\vee$ and $\times$ realizes $\cup : \mathcal{A}(\mathbf{X}) \times \mathcal{A}(\mathbf{X}) \to \mathcal{A}(\mathbf{X})$.
\item This follows from 1. and 2. using de Morgan's law.
\item Composition of $\bigvee$ with the input yields the output, and is computable due to Proposition \ref{compana:prop:functionspacesbasics} (4).
\item This follows from 1. and 4. using de Morgan's law.
\item Again, this is a special case of composition, which is computable due to Proposition \ref{compana:prop:functionspacesbasics} (4).
\item Here we have a special case of $\eval$, which is computable due to Proposition \ref{compana:prop:functionspacesbasics} (1).
\item This follows from composing the computable function $\times : \mathcal{C}(\mathbf{X}, \mathbb{S}) \times \mathcal{C}(\mathbf{Y}, \mathbb{S}) \to \mathcal{C}(\mathbf{X} \times \mathbf{Y}, \mathbb{S} \times \mathbb{S})$ from Proposition \ref{compana:prop:functionspacesbasics} (6) together with $\wedge : \mathbb{S} \times \mathbb{S} \to \mathbb{S}$ and $\vee : \mathbb{S} \times \mathbb{S} \to \mathbb{S}$ respectively.
\item This follows from Proposition \ref{compana:prop:functionspacesbasics} (1,2) via $x \in \operatorname{Cut}(y, U)$ iff $(x, y) \in U$.
\item Essentially, this is composition with $\bigvee : \mathcal{C}(\mathbb{N}, \mathbb{S}) \to \mathbb{S}$ from Proposition \ref{prop:sierpbasics}.
\end{enumerate}
\end{proof}
\end{proposition}

A represented space $\mathbf{X}$ canonically induces a topological space by equipping the set $X$ by the quotient topology $\mathcal{T}_\mathbf{X}$ of $\delta_X : \subseteq \Cantor \to X$. One can verify that $\mathcal{T}_\mathbf{X}$ is the underlying set of $\mathcal{O}(\mathbf{X})$; i.e.~that the open sets of a represented space are the open sets in the induced topological space: Essentially, the open subset of a subspace of $\Cantor$ that witnesses that a set $U$ is in $\mathcal{T}_\mathbf{X}$ can be translated into a continuous realizer of $\chi_U : \mathbf{X} \to \mathbb{S}$. Proposition \ref{compana:prop:basicsetoperations} (6) implies that a continuous map between represented spaces is also continuous as a map between the induced topological spaces. However, the converse is generally false.

\subsubsection*{Remarks}
The introduction of \Sierp space to the study of representations, and subsequent use to derive representations of the open and the closed subsets, is due to \name{Schr\"oder} in 2002 \cite{schroder5} (cf.~\cite{presser}). However, the usefulness of such definitions to obtain concise proofs has not been fully appreciated in this setting (e.g.~\name{Weihrauch} and \name{Grubba} are not using it in \cite{grubba} from 2009). Statements equivalent to the items of Proposition \ref{compana:prop:basicsetoperations} have been proven in various restricted settings before.

Using the combination of the function space construction and the presence of \Sierp space to obtain the space of open (closed) sets together with the standard operations on them is the fundamental idea of synthetic topology (\name{Escard\'o} 2004 \cite{escardo}).

\name{Collins} has stated the results of Proposition \ref{compana:prop:basicsetoperations} $(1.-6.,9.)$ as part of \cite[Theorem 3.23, Theorem 3.24, \& Proposition 3.26]{collins5}.

The space $\mathcal{A}(\mathbf{X})$ has been studied as the upper Fell topology on the hyperspace of closed sets in topology. This topology was introduced by \name{Fell} in 1962 \cite{fell}. Note that the characterization of compact sets in Section \ref{sec:compactness} is crucial for this coincidence, which has been observed before by \name{Brattka} and \name{Presser} 2003 \cite{presser}. A general source for hyperspace topologies is \cite{beer}.
\section{Compactness}
\label{sec:compactness}
Compactness is occasionally described as a generalization of finiteness, in that compactness means that it can be verified that a property applies to all elements of a space -- so in same sense, compact spaces admit a form of exhaustive search\footnote{For the idea of exhaustive search of infinite sets, see \name{Escard\'o} \cite{escardo5}.}, despite potentially having uncountable cardinality. Compactness both plays a r\^ole as a property of spaces, as well as inducing a represented space of subsets of a given space. We provide various equivalent characterizations of the former, and list various useful computable operations on the latter.
\begin{definition}
\label{compana:def:compactness}
A represented space $\mathbf{X}$ is (computably) compact, if the map $\isempty_{\mathbf{X}} : \mathcal{A}(\mathbf{X}) \to \mathbb{S}$ defined by $\isempty_{\mathbf{X}}(\emptyset) = \top$ and $\isempty_{\mathbf{X}}(A) = \bot$ otherwise is continuous (computable).
\end{definition}

\begin{proposition}
\label{compana:prop:compactnessbasics}
The following properties are equivalent for a represented space $\mathbf{X}$:
\begin{enumerate}
\item $\mathbf{X}$ is (computably) compact.
\item $\operatorname{IsFull}_{\mathbf{X}} : \mathcal{O}(\mathbf{X}) \to \mathbb{S}$ defined by $\operatorname{IsFull}_{\mathbf{X}}(X) = \top$ and $\operatorname{IsFull}_{\mathbf{X}}(U) = \bot$ otherwise is continuous (computable).
\item For every (computable\footnote{We remind the reader that we call $A \in \mathcal{A}(\mathbf{X})$ computable, if it is a computable name. This generalizes Weihrauch's notion of a co-c.e. closed set, not that of a computable closed set!}) $A \in \mathcal{A}(\mathbf{X})$ the subspace $\mathbf{A}$ is (computably) compact.
\item $\subseteq : \mathcal{A}(\mathbf{X}) \times \mathcal{O}(\mathbf{X}) \to \mathbb{S}$ defined by $\subseteq(A, U) = \top$, iff $A \subseteq U$ is continuous (computable).
\item $\operatorname{IsCover} : \mathcal{C}(\mathbb{N}, \mathcal{O}(\mathbf{X})) \to \mathbb{S}$ defined by $\operatorname{IsCover}((U_n)_{n \in \mathbb{N}}) = \top$, iff $\bigcup_{n \in \mathbb{N}} U_n = X$ is continuous (computable).
\item $\operatorname{FiniteSubcover} : \subseteq \mathcal{C}(\mathbb{N}, \mathcal{O}(\mathbf{X})) \mto \mathbb{N}$ with $\dom(\operatorname{FiniteSubcover}) \linebreak = \{(U_n)_{n \in \mathbb{N}} \mid \exists N \bigcup_{n \leq N} U_n = X\}$ and $N \in \operatorname{FiniteSubcover}((U_n)_{n \in \mathbb{N}})$ iff $\bigcup_{n \leq N} U_n = X$ is continuous (computable).
\item $\operatorname{Enough} :\subseteq \mathcal{C}(\mathbb{N}, \mathcal{A}(\mathbf{X})) \mto \mathbb{N}$ with $(A_i)_{i \in \mathbb{N}} \in \dom(\operatorname{Enough})$ iff $\exists N \ \bigcap_{i \leq N} A_i = \emptyset$, and $N \in \operatorname{Enough}((A_i)_{i \in \mathbb{N}})$ iff $\bigcap_{i \leq N} A_i = \emptyset$ is continuous (computable).
\item For every represented spaces $\mathbf{Y}$, the map $\pi_2 : \mathcal{A}(\mathbf{X} \times \mathbf{Y}) \to \mathcal{A}(\mathbf{Y})$ defined by $\pi_2(A) = \{y \in \mathbf{Y} \mid \exists x \in \mathbf{X} \ (x, y) \in A\}$ is well-defined and continuous (computable).
\item For some non-empty represented space $\mathbf{Y}$ (containing a computable point), the map $\pi_2 : \mathcal{A}(\mathbf{X} \times \mathbf{Y}) \to \mathcal{A}(\mathbf{Y})$ is well-defined and continuous (computable).
\end{enumerate}
\begin{proof}
\begin{description}
\item[$1. \Leftrightarrow 2.$] $\isempty_{\mathbf{X}}$ and $\operatorname{IsFull}_{\mathbf{X}}$ have exactly the same realizers.
\item[$1. \Leftrightarrow 3.$] $X$ is always a computable element of $\mathcal{A}(\mathbf{X})$, as it is realized by the constant function $p \mapsto 0^\mathbb{N}$. This provides the implication $1. \Leftarrow 3.$. For the other direction, note $\isempty_{\mathbf{A}}(B) = \isempty_{\mathbf{X}}(A \cap B)$, so continuity (computability) of $\isempty_{\mathbf{A}}$ follows from Proposition \ref{compana:prop:basicsetoperations} (2) together with Corollary \ref{compana:corr:partial} and Proposition \ref{compana:prop:functionspacesbasics} (2).
\item[$1. \Leftrightarrow 4.$] $\emptyset \in \mathcal{O}(\mathbf{X})$ is computable by Proposition \ref{compana:prop:functionspacesbasics} (6), so if $\subseteq$ is computable, so is $A \mapsto \subseteq(A, \emptyset)$ by Corollary \ref{compana:corr:partial}. For the other direction, observe $A \subseteq B \Leftrightarrow A \cap B^C = \emptyset$, so we may use the combination of Proposition \ref{compana:prop:basicsetoperations} (1, 3), Corollary \ref{compana:corr:partial} and Proposition \ref{compana:prop:functionspacesbasics} (2) to obtain computability of $\subseteq$ here.
\item[$2. \Leftrightarrow 5.$] By Proposition \ref{compana:prop:functionspacesbasics} (6) we can compute $(U)_{n \in \mathbb{N}}$ from $U$, and clearly $\operatorname{IsFull}(U) = \operatorname{IsCover}((U)_{n \in \mathbb{N}})$. On the other hand, by Proposition \ref{compana:prop:basicsetoperations} (4), we can compute $\bigcup_{n \in \mathbb{N}} U_n$ from $(U_n)_{n \in \mathbb{N}}$, and $\operatorname{IsCover}((U)_{n \in \mathbb{N}}) = \operatorname{IsFull}(\bigcup_{n \in \mathbb{N}} U_n)$.
\item[$5. \Rightarrow 6.$] Assume that the Turing machine $M$ computes $\operatorname{IsCover}$ (potentially with access to some oracle). In order to solve $\operatorname{FiniteSubcover}$, we simulate $M$ on the input for $\operatorname{FiniteSubcover}$ (which is of a suitable type). As we know that the input sequence \emph{does} cover $X$, we also know that $M$ has to write a $1$ eventually. When $M$ writes the first $1$, it has only read some finite prefix of the input. In particular, we may assume that $M$ has no information at all about the $U_n$ with $n > N$ for some $N \in \mathbb{N}$. Moreover, we can find such an $N$ effectively from observing the simulation of $M$.

    Now $N$ constitutes a valid answer to $\operatorname{FiniteSubcover}((U_n)_{n \in \mathbb{N}})$. To see this, assume the contrary. Then $\bigcup_{n = 0}^N U_n \neq X$. Now consider the sequence $(U'_n)_{n \in \mathbb{N}}$ with $U'_n = U_n$ for $n \leq N$ and $U'_n = \emptyset$ otherwise. On some name for $(U'_n)_{n \in \mathbb{N}}$, $M$ will eventually print a $1$ - as it cannot distinguish $(U'_n)_{n \in \mathbb{N}}$ from $(U_n)_{n \in \mathbb{N}}$ before that. But as $(U'_n)_{n \in \mathbb{N}}$ does not constitute a cover of $X$, this contradicts the initial assumption.
\item[$6. \Rightarrow 2.$] Assume that the Turing machine $M$ computes $\operatorname{FiniteSubcover}$ (potentially with access to some oracle). In order to solve $\operatorname{IsFull}(U)$, we simulate $M$ on input $(U)_{n \in \mathbb{N}}$, which we can obtain by Proposition \ref{compana:prop:functionspacesbasics} (6). Beside the simulation, we repeatedly write $0$s on the output tape. If $M$ ever produces some $N \in \mathbb{N}$ as output, we write a $1$. This actually solves $\operatorname{IsFull}(U)$ correctly.

    If $U = X$, then $(U)_{n \in \mathbb{N}}$ is a valid input for $\operatorname{FiniteSubcover}$, so eventually some $N \in \mathbb{N}$ will be produced, causing the final result to be $\top \in \mathbb{S}$. Now assume that $M$ produces some $N \in \mathbb{N}$ on input $(U)_{n \in \mathbb{N}}$ for $U \neq X$. At the time where $N$ has been written, $M$ has only received information about some finite prefix of its input. In particular, there is some $K > N$ s.t.~ $M$ exhibits the same behaviour when faced with the input sequence $(U'_n)_{n \in \mathbb{N}}$ with $U'_n = U$ for $n \leq K$ and $U'_n = X$ otherwise. As $(U'_n)_{n \in \mathbb{N}}$  is a valid input for $\operatorname{FiniteSubcover}$, $M$ would be required to produce some $k \in \mathbb{N}$ with $k \geq K$ rather than $N$. The contradiction can only be resolved by the assumption that $M$ never completes an output on input $(U)_{n \in \mathbb{N}}$ for $U \neq X$, which yields the final answer to be $\bot \in \mathbb{S}$.
\item[$6. \Leftrightarrow 7.$] This follows via Proposition \ref{compana:prop:basicsetoperations} (1).
\item[$1. \Rightarrow 8.$] By Proposition \ref{compana:prop:functionspacesbasics} (2), we can compute $y \mapsto (x \mapsto A(x, y))$ from $A \in \mathcal{A}(\mathbf{X} \times \mathbf{Y})$. Now we find $(x \mapsto A(x, y)) \in \mathcal{A}(\mathbf{X})$, and if $\isempty_\mathbf{X}$ is computable, so is $y \mapsto \isempty_\mathbf{X}(x \mapsto A(x, y)) = \pi_2(A)$.
\item[$8. \Rightarrow 9.$] Just instantiate with $\mathbf{Y} = \mathbb{N}$.
\item[$9. \Rightarrow 1.$] Let $y \in \mathbf{Y}$ be a (computable) point. Then $A \mapsto (y \in \pi_2(A \times Y)^C)$ realizes $\isempty_\mathbf{X}$, and is continuous (computable) using the assumption, together with Proposition \ref{compana:prop:basicsetoperations} (1, 6, 7) together with Corollary \ref{compana:corr:partial}.
\end{description}
\end{proof}
\end{proposition}

\begin{proposition}
\label{compana:prop:compactnessimage1}
Let $\mathbf{X}$ be (computably) compact and $f : \mathbf{X} \to \mathbf{Y}$ a (computable) continuous surjection. Then $\mathbf{Y}$ is (computably) compact.
\begin{proof}
We use the characterization of compactness provided by Proposition \ref{compana:prop:compactnessbasics} (2). By Proposition \ref{compana:prop:basicsetoperations} (6) we find $f^{-1} : \mathcal{O}(\mathbf{Y}) \to \mathcal{O}(\mathbf{X})$ to be continuous (computable). Now observe $\operatorname{IsFull}_{\mathbf{Y}}(U) = \operatorname{IsFull}_{\mathbf{X}}(f^{-1}(U))$ due to surjectivity of $f$.
\end{proof}
\end{proposition}

\begin{proposition}
\label{prop:compactproduct}
If $\mathbf{X}$, $\mathbf{Y}$ are (computably) compact, then so is $\mathbf{X} \times \mathbf{Y}$.
\begin{proof}
Note $\isempty_{\mathbf{X} \times \mathbf{Y}} = \isempty_\mathbf{Y} \circ \pi_2$, and use Proposition \ref{compana:prop:compactnessbasics} (8).
\end{proof}
\end{proposition}

We can call a subset of a represented space compact, if it is compact when represented with the subspace representation. A set $K \subseteq \mathbf{X}$ is called \emph{saturated}, iff $K = \bigcap_{\{U \in \mathcal{O}(\mathbf{X}) \mid K \subseteq U\}} U$. Noting that for compact $K$ the set $\{U \in \mathcal{O}(\mathbf{X}) \mid K \subseteq U\}$ is open in $\mathcal{O}(\mathbf{X})$, we obtain a representation of the set $\mathcal{K}(\mathbf{X})$ of saturated compact sets by identifying it as a subspace of $\mathcal{O}(\mathcal{O}(\mathbf{X}))$. In particular, this makes $\operatorname{IsContainedIn}  : \mathcal{K}(\mathbf{X}) \times \mathcal{O}(\mathbf{X}) \to \mathbb{S}$ computable, and provides a uniform counterpart to the results in \ref{compana:prop:compactnessbasics}. For arbitrary sets $K \subseteq X$, we use $\sat{K} := \bigcap_{\{U \in \mathcal{O}(\mathbf{X}) \mid K \subseteq U\}} U$ to denote its \emph{saturation}, and point out that any $\sat{K}$ is saturated. In a $T_1$ space, i.e.~a represented space $\mathbf{X}$ where $\forall x \in \mathbf{X} \ \{x\} \in \mathcal{A}(\mathbf{X})$, any set is already saturated.

We proceed to exhibit a number of computable operations on spaces of (saturated) compact sets, some of which can be seen as uniform counterparts of results about compact spaces above.
\begin{proposition}
\label{prop:operationsoncompacts}
The following operations are well-defined and computable:
\begin{enumerate}
\item $\operatorname{IsContainedIn} : \mathcal{K}(\mathbf{X}) \times \mathcal{O}(\mathbf{X}) \to \mathbb{S}$
\item $x \mapsto \sat{\{x\}} : \mathbf{X} \to \mathcal{K}(\mathbf{X})$
\item $\cup : \mathcal{K}(\mathbf{X}) \times \mathcal{K}(\mathbf{X}) \to \mathcal{K}(\mathbf{X})$
\item $\sat{\cap} : \mathcal{K}(\mathbf{X}) \times \mathcal{A}(\mathbf{X}) \to \mathcal{K}(\mathbf{X})$
\item $f^{-1} \mapsto \sat{f} : \subseteq \mathcal{C}(\mathcal{O}(\mathbf{Y}), \mathcal{O}(\mathbf{X})) \to \mathcal{C}(\mathcal{K}(\mathbf{X}), \mathcal{K}(\mathbf{Y}))$; here and below $\sat{f}$ is the composition of $f$ and the saturation operator $\sat$
\item $f \mapsto \sat{f} : \mathcal{C}(\mathbf{X}, \mathbf{Y}) \to \mathcal{C}(\mathcal{K}(\mathbf{X}), \mathcal{K}(\mathbf{Y}))$
\item $(f, K) \mapsto \sat{f[K]} : \mathcal{C}(\mathbf{X}, \mathbf{Y}) \times \mathcal{K}(\mathbf{X}) \to \mathcal{K}(\mathbf{Y})$
\item $\times : \mathcal{K}(\mathbf{X}) \times \mathcal{K}(\mathbf{Y}) \to \mathcal{K}(\mathbf{X} \times \mathbf{Y})$
\item $\pi_1 : \mathcal{K}(\mathbf{X} \times \mathbf{Y}) \to \mathcal{K}(\mathbf{X})$, $\pi_2 : \mathcal{K}(\mathbf{X} \times \mathbf{Y}) \to \mathcal{K}(\mathbf{Y})$
\end{enumerate}
\begin{proof}
\begin{enumerate}
\item Taking into account the definition of $\mathcal{K}$, this is an instantiation of $\eval$ from Proposition \ref{compana:prop:functionspacesbasics} (1).
\item Note that $x \in U$ iff $\sat{\{x\}} \subseteq U$ for $U \in \mathcal{O}(\mathbf{X})$.
\item $\cup : \mathcal{K}(\mathbf{X}) \times \mathcal{K}(\mathbf{X}) \to \mathcal{K}(\mathbf{X})$ is realized by $\cap : \mathcal{O}(\mathcal{O}(\mathbf{X})) \times \mathcal{O}(\mathcal{O}(\mathbf{X})) \to \mathcal{O}(\mathcal{O}(\mathbf{X}))$ from Proposition \ref{compana:prop:basicsetoperations} (2).
\item The core observation is that $A \cap B \subseteq U$ iff $A \subseteq (U \cup B^C)$ together with Proposition \ref{compana:prop:basicsetoperations} (1, 2).
\item It suffices to show that $f^{-1} \mapsto f$ is well-defined, its computability then follows from the computability of $f^{-1} \mapsto (f^{-1})^{-1} : \subseteq \mathcal{C}(\mathcal{O}(\mathbf{Y}), \mathcal{O}(\mathbf{X})) \to \mathcal{C}(\mathcal{O}(\mathcal{O}(\mathbf{X})), \mathcal{O}(\mathcal{O}(\mathbf{Y})))$ as it is a restriction of this map. The computability of the latter map in turn is a special case of Proposition \ref{compana:prop:basicsetoperations} (6). Well-definedness in turn follows directly from the characterization of $\mathcal{K}(\mathbf{X})$ as subspace of $\mathcal{O}(\mathcal{O}(\mathbf{X}))$.
\item From (5) together with Proposition \ref{compana:prop:basicsetoperations} (6).
\item From (6) via type conversion (Proposition \ref{compana:prop:functionspacesbasics}).
\item We show that given $A \in \mathcal{K}(\mathbf{X})$, $B \in \mathcal{K}(\mathbf{Y})$ and $U \in \mathcal{O}(\mathbf{X} \times \mathbf{Y})$, we can semidecide $A \times B \subseteq U$. Drawing from (1) and Proposition \ref{compana:prop:basicsetoperations} (9), we consider $y \mapsto \operatorname{IsContainedIn}(A, \operatorname{Cut}(y, U))$, which defines some element $U_A$ of $\mathcal{O}(\mathbf{Y})$. Furthermore, we note that $A \times B \subseteq U$ iff $B \subseteq U_A$.
\item By Proposition \ref{prop:categorystructure} the maps $\pi_1 : \mathbf{X} \times \mathbf{Y} \to \mathbf{X}$, $\pi_2 : \mathbf{X} \times \mathbf{Y} \to \mathbf{Y}$ are computable. Corollary \ref{compana:corr:partial} allows us to use (6) to lift this to compact sets.
\end{enumerate}
\end{proof}
\end{proposition}

\begin{corollary}
\label{compana:corr:idclosedcompact}
$\mathbf{X}$ is (computably) compact, iff $\sat{\id} : \mathcal{A}(\mathbf{X}) \to \mathcal{K}(\mathbf{X})$ is well-defined and continuous (computable).
\end{corollary}

\subsection*{Remarks}
The represented space $\mathcal{K}(\mathbf{X})$ was introduced, and various parts of Proposition \ref{prop:operationsoncompacts} were proven by Schr\"oder in \cite[Section 4.4.3]{schroder5} 2002 (for admissible $\mathbf{X}$). This definition of compactness also follows the approach of synthetic topology (\name{Escard\'o} 2004 \cite{escardo}), and some of our results have also been obtained there (Proposition \ref{compana:prop:compactnessbasics} $[2. \Rightarrow 3.]$, Propositions \ref{compana:prop:compactnessimage1}, \ref{prop:compactproduct}). The equivalence in Proposition \ref{compana:prop:compactnessbasics} $[2. \Leftrightarrow 8.]$ was shown by \name{Escard\'o} in \cite{escardo2}.

That $\mathcal{K}(\mathbf{X})$ can be identified with a certain subspace of $\mathcal{O}(\mathcal{O}(\mathbf{X}))$ is the statement of the Hofman-Mislove theorem (e.g.~\cite{mislove}, see also the extension by \name{Schr\"oder} in \cite{schroder7}).

The role of saturation as canonization operation for compact sets is already discussed by \name{Collins} \cite{collins5}. \name{Collins} has stated the results of Proposition \ref{prop:operationsoncompacts} $(2.-4.,6.,7.)$ as part of \cite[Theorem 3.23 \& Theorem 3.24]{collins5}; Corollary \ref{compana:corr:idclosedcompact} as \cite[Theorem 3.31 (4)]{collins5}; and Proposition \ref{compana:prop:compactnessbasics} ($1. \Rightarrow 8.$) as \cite[Proposition 3.33]{collins5}. His definition of compactness \cite[Definition 3.28 (4)]{collins5} is equivalent to the one here. However, the claimed characterization of the computably compact spaces as the images of Cantor space computable functions in \cite[Proposition 3.30]{collins5} (stated without proof) is wrong. A counterexample (due to \name{de Brecht}, personal communication) is one-point compactification of $\Baire$ (\footnote{In general, the one-point compactification of a represented space $(X, \delta)$ can be introduced as $(X \cup \{\bot\}, \delta_C)$ where $\delta_C(0^n1p) = \delta(p)$ and $\delta_C(0^\mathbb{N}) = \bot$. Any one-point compactification is computably compact, but can only have a total Cantor space representation if the original space was locally compact.}). This also impacts \cite[Theorem 3.32]{collins5}.

Some results in this section generalize known results in the far more restricted setting of computable metric spaces. A version of Proposition \ref{compana:prop:compactnessbasics} $[5. \Leftrightarrow 6.]$ was proven by \name{Brattka} and \name{Presser} (2003 \cite{presser}). \name{Weihrauch} proved the restricted version of \ref{prop:operationsoncompacts} (7) in 2003 \cite{weihrauchf}. For computable metric spaces, compact compactness can also be characterized by the computability of the outer radius of closed sets (taking values in $\mathbb{R}_>$) as shown by \name{Le Roux} and the author as \cite[Proposition 19]{paulyleroux-arxiv}.

The equivalences in Proposition \ref{compana:prop:compactnessbasics} $[3. \Leftrightarrow 6. \Leftrightarrow 7. \Leftrightarrow 8.]$ are uniform counterparts to well-known characterizations of compactness in topology.

The space $\mathcal{K}(\mathbf{X})$ corresponds to the upper Vietoris topology (which is often defined on the hyperspace of closed sets). This topology was introduced by \name{Vietoris} in 1922 \cite{vietoris}. Again, a general source for hyperspace topologies is \cite{beer}.

\section{$T_2$ separation}
The $T_2$ separation axiom can, in the context of represented spaces, be understood equivalently as the property of a space making either inequality or the subspace of compact sets well-behaved. The strong connection between $T_2$ separation and compactness is somewhat reminiscent of the ultrafilter approach to topology, where compactness means each ultrafilter converges to at least one point, and $T_2$ that they converge to at most one point. Several of our equivalences will require the $T_0$ property to work, which we understand in a non-effective way to mean that $x \neq y$ implies $\sat{\{x\}} \neq \sat{\{y\}}$. We also need the (non-effective and non-uniform) $T_1$ property, which we understand to mean that if $x \neq y$, then $\sat{\{x\}} \cap \sat{\{y\}} = \emptyset$.

\begin{definition}
\label{compana:def:t1}
A represented space $\mathbf{X}$ is (computably) $T_2$, if the map $x \mapsto \{x\} : \mathbf{X} \to \mathcal{A}(\mathbf{X})$ is well-defined and continuous (computable).
\end{definition}

\begin{proposition}
\label{compana:prop:t1characterization}|
The following properties are equivalent for a represented space $\mathbf{X}$:
\begin{enumerate}
\item $\mathbf{X}$ is (computably) $T_2$.
\item $\id : \mathcal{K}(\mathbf{X}) \to \mathcal{A}(\mathbf{X})$ is well-defined and continuous (computable), and $\mathbf{X}$ is $T_0$.
\item $\cap : \mathcal{K}(\mathbf{X}) \times \mathcal{K}(\mathbf{X}) \to \mathcal{K}(\mathbf{X})$ is well-defined and continuous (computable), and $\mathbf{X}$ is $T_1$.
\item $x \mapsto \sat{\{x\}} : \mathbf{X} \to \mathcal{A}(\mathbf{X})$ is well-defined, injective and continuous (computable).
\item $\mathalpha{\neq} : \mathbf{X} \times \mathbf{X} \to \mathbb{S}$ defined by $\mathalpha{\neq}(x, x) = \bot$ and $\mathalpha{\neq}(x, y) = \top$ otherwise is continuous (computable).
\item $\Delta_\mathbf{X} = \{(x, x) \mid x \in \mathbf{X}\} \in \mathcal{A}(\mathbf{X} \times \mathbf{X})$ (is computable\footnote{Again, we remind the reader that we call $A \in \mathcal{A}(\mathbf{X})$ computable, if it is a computable name. This generalizes Weihrauch's notion of a co-c.e. closed set, not that of a computable closed set!}).
\item $\operatorname{Graph} : \mathcal{C}(\mathbf{Y}, \mathbf{X}) \to \mathcal{A}(\mathbf{Y} \times \mathbf{X})$ is well-defined and continuous (computable) for any represented space $\mathbf{Y}$.
\end{enumerate}
\begin{proof}
\begin{description}
\item[$1. \Rightarrow 2.$] The map is identical to $K \mapsto (x \mapsto \operatorname{IsContainedIn}((K \cap \{x\}), \emptyset))$. To see that the map is continuous (computable) under the assumption that $x \mapsto \{x\} : \mathbf{X} \to \mathcal{A}(\mathbf{X})$ is continuous (computable), first use Propositions \ref{prop:operationsoncompacts} (4) to obtain $K \cap \{x\} \in \mathcal{K}(\mathbf{X})$, then Proposition\ref{prop:operationsoncompacts} (1). Being $T_0$ is an obvious consequence.
\item[$1. \wedge 2. \Rightarrow 3.$] The continuity (computability) of the map in $3.$ is a direct consequence of the continuity (computability) of the map in $2.$ and Proposition  \ref{prop:operationsoncompacts} (4). The $T_1$ property follows from the well-definedness of $x \mapsto \{x\} : \mathbf{X} \to \mathcal{A}(\mathbf{X})$.
\item[$3. \Rightarrow 2.$] Observe that $x \in K$ holds for a compact set $K$, if and only if $\sat{\{x\}} \subseteq K$ holds. The $T_1$-property allows us to strengthen this to $x \in K$ iff $\sat{\{x\}} \cap K \neq \emptyset$. Hence $\id : \mathcal{K}(\mathbf{X}) \to \mathcal{A}(\mathbf{X})$ is given by $K \mapsto (x \mapsto \operatorname{IsContainedIn}(K \cap \sat{\{x\}}, \emptyset))$, if intersection of compact sets is continuous (computable) (using Propositions \ref{prop:operationsoncompacts} (1)).
\item[$2. \Rightarrow 4.$] Continuity (computability) and well-definedness of the map in $4.$ follows immediately from continuity (computability) and well-definedness of the map in $2.$ and Proposition  \ref{prop:operationsoncompacts} (2). Injectiveness of the map is equivalent to the $T_0$-property.
\item[$4. \Rightarrow 2.$] Similar to $3. \Rightarrow 2.$: In $K \mapsto (x \mapsto \operatorname{IsContainedIn}(K \cap \sat{\{x\}}, \emptyset))$ use $x \mapsto \sat{\{x\}} : \mathbf{X} \to \mathcal{A}(\mathbf{X})$, and the intersection from Proposition \ref{prop:operationsoncompacts} (4).
\item[$2. \Rightarrow 5.$] Given $x, y \in \mathbf{X}$, compute $\sat{\{x\}}, \sat{\{y\}} \in \mathcal{K}(\mathbf{X})$ by Proposition \ref{prop:operationsoncompacts} (2), and then $\sat{\{x\}}, \sat{\{y\}} \in \mathcal{A}(\mathbf{X})$ by the assumption. Now the claim follows from $x \neq y$ iff $x \notin \sat{\{y\}} \vee y \notin \sat{\{x\}}$. The equivalence holds due to the $T_0$ property, and the right hand side is computable due to Proposition \ref{prop:sierpbasics}.
\item[$5. \Rightarrow 1.$] By Proposition \ref{compana:prop:functionspacesbasics} (2) we find $x \mapsto (y \mapsto \mathalpha{\neq}(x, y))$ to be continuous (computable), but this has the same realizers as $x \mapsto \{x\}$.
\item[$5. \Leftrightarrow 6.$] This is just a reformulation along the definition of $\mathcal{A}(\mathbf{X} \times \mathbf{X})$.
\item[$5. \Rightarrow 7.$] This follows from $\operatorname{Graph}(f) = \{(y, x) \mid f(y) = x\} \in \mathcal{A}(\mathbf{Y} \times\mathbf{X})$; that $f \times \id : \mathbf{Y} \times \mathbf{X} \to \mathbf{X} \times \mathbf{X}$ is available as a continuous function; and Proposition \ref{compana:prop:basicsetoperations} (1,6).
\item[$7. \Rightarrow 6.$] Pick $\mathbf{Y} := \mathbf{X}$, and then consider $\operatorname{Graph}(\id_\mathbf{X})$.
\end{description}
\end{proof}
\end{proposition}

At the first glance, Definition \ref{compana:def:t1} seems to be a prime candidate for $T_1$ separation rather than $T_2$ separation. In light of Proposition \ref{compana:prop:t1characterization}, it is clear that the denotation $T_2$ is justified, too. The reason to prefer the identification as $T_2$ rather than $T_1$ separation lies in the observation that a represented space is $T_2$, if and only if the induced topological space is sequentially $T_2$ (a topological space is sequentially $T_2$ iff its diagonal is sequentially closed (e.g.~\cite[Section 2.4]{schroder7}), so the claim is Proposition \ref{compana:prop:t1characterization} (6)). On the other hand, there are represented spaces not satisfying Definition \ref{compana:def:t1} but inducing a $T_1$ topology\footnote{An example is the subspace $\{\{n\} \mid n \in \mathbb{N}\} \subseteq \mathcal{A}(\mathbb{N})$, as pointed out by \name{Weihrauch} in \cite{weihrauchh}. This space corresponds to the cofinite topology on $\mathbb{N}$.}.

With both compactness and the computable $T_2$ property in place, we can proceed to discuss \emph{computably proper maps}, i.e. those (continuous) maps where the preimages of compacts are compact. First, we present a counterpart to a classic result from topology:

\begin{proposition}
Let $\mathbf{X}$ be (computably) compact and $\mathbf{Y}$ (computably) $T_2$. Then the continuous maps in $\mathcal{C}(\mathbf{X}, \mathbf{Y})$ are uniformly proper, i.e. the map $(f, K) \mapsto f^{-1}(K) : \mathcal{C}(\mathbf{X}, \mathbf{Y}) \times \mathcal{K}(\mathbf{Y}) \to \mathcal{K}(\mathbf{X})$ is well-defined and continuous (computable).
\begin{proof}
By Proposition \ref{compana:prop:t1characterization} (2) the map $\id : \mathcal{K}(\mathbf{Y}) \to \mathcal{A}(\mathbf{Y})$ is continuous (computable) due to the $T_2$-property for $\mathbf{Y}$. The map $(f, A) \mapsto f^{-1}(A) : \mathcal{C}(\mathbf{X}, \mathbf{Y}) \times \mathcal{A}(\mathbf{Y}) \to \mathcal{A}(\mathbf{X})$ is computable by Proposition \ref{compana:prop:basicsetoperations} (1, 6), Proposition \ref{compana:prop:functionspacesbasics} (1). By Corollary \ref{compana:corr:idclosedcompact} the map $\id : \mathcal{A}(\mathbf{X}) \to \mathcal{K}(\mathbf{X})$ is continuous (computable) due to the compactness of $\mathbf{X}$. Composition of these maps yields the claim.
\end{proof}
\end{proposition}

\begin{proposition}
Let $f : \mathbf{X} \to \mathbf{Y}$ be surjective, continuous (computable) and (computably) proper, i.e. let $f^{-1} : \mathcal{K}(\mathbf{Y}) \to \mathcal{K}(\mathbf{X})$ be well-defined and continuous (computable). If $\mathbf{X}$ is (computably) $T_2$ and $\mathbf{Y}$ is $T_1$, then $\mathbf{Y}$ even is (computably) $T_2$.
\begin{proof}
We use the characterization of the $T_2$ property given in Proposition \ref{compana:prop:t1characterization} (3), i.e. the intersection of compact sets being computable as a compact set. For surjective $f$, we have $A \cap B = f[f^{-1}(A) \cap f^{-1}(B)]$. Using this equation for $A, B \in \mathcal{K}(\mathbf{Y})$, we can compute $f^{-1}(A)$, $f^{-1}(B)$ by the assumption $f$ were computably proper, then $f^{-1}(A) \cap f^{-1}(B)$ as a compact set by assumption $\mathbf{X}$ were $T_2$, and finally $\sat{f[f^{-1}(A) \cap f^{-1}(B)]} \in \mathcal{K}(\mathbf{Y})$ due to Proposition \ref{prop:operationsoncompacts} (6). Thus, we see that we can obtain $\sat{(A \cap B)} \in \mathcal{K}(\mathbf{Y})$ from $A, B \in \mathcal{K}(\mathbf{Y})$. By inspecting the proof of Proposition \ref{compana:prop:t1characterization} ($3. \Rightarrow 2.$), we notice that we only care whether the intersection of two compact sets is empty or not. As saturation preserves the empty set, we conclude that computability of $\sat{\cap} : \mathcal{K}(\mathbf{Y}) \times \mathcal{K}(\mathbf{Y}) \to \mathcal{K}(\mathbf{Y})$ suffices for a $T_1$ space $\mathbf{Y}$ to be computably Hausdorff.
\end{proof}
\end{proposition}

\begin{corollary}
\label{corr:propert2}
If $\mathbf{X}$ is $T_1$ and has a computably proper representation $\delta_\mathbf{X}$, then $\mathbf{X}$ is computably $T_2$.
\end{corollary}

\subsubsection*{Remarks}
\name{Escard\'o} uses the condition Proposition \ref{compana:prop:t1characterization} (5) to introduce the Hausdorff property in \cite{escardo}; and proceeds to prove Proposition \ref{compana:prop:t1characterization} $[5. \Rightarrow 2.]$. \name{Weihrauch} studied computable separation in countably based spaces in \cite{weihrauchh} (2010), including what corresponds to the equivalence Proposition \ref{compana:prop:t1characterization} $[1. \Rightarrow 5.]$.

\name{Collins} translated \name{Escard\'o}'s definition of Hausdorff and the equivalence from Proposition \ref{compana:prop:t1characterization} $[5. \Rightarrow 2.]$ into computable analysis in \cite[Definition 3.28 (2 \& Theorem 3.31 (2.)]{collins5} (under the name \emph{effectively distinguishable}).

Proper representations were studied by \name{Schr\"oder} in \cite{schroder6} 2004; in particular, \cite[Theorem 5.3]{schroder6} shows that having a proper admissible representation even implies metrizability (and is implied by metrizability, in turn); which is much stronger than Corollary \ref{corr:propert2}. A minor further strengthening of the result can be found in \cite{pauly-kihara-arxiv}, where it is shown that any such space embeds computability into a computable metric space.

\section{Overtness}
\label{sec:overt}
Like compactness, \emph{overtness} both is a property of represented spaces and induces a space of certain subsets. Unlike compactness, overtness is classically vacuous and thus not that well-known. As discussed below, overtness has been present in the study of representations of spaces of subsets in computable analysis for a long time, although in a disguised r\^ole.
\begin{definition}
\label{def:overtspace}
Let $\nonempty_{\mathbf{X}} : \mathcal{O}(\mathbf{X}) \to \mathbb{S}$ be defined by $\nonempty_\mathbf{X}(\emptyset) = \bot$ and $\nonempty_\mathbf{X}(U) = \top$ for $U \neq \emptyset$. Now call $\mathbf{X}$ (computably) \emph{overt}, iff $\nonempty_{\mathbf{X}}$ is continuous (computable).
\end{definition}

This definition shows that overtness is in some sense the dual to compactness. In computable analysis, typically the criterion separable rather than overt has been used for spaces, based on the following observation:
\begin{proposition}
\label{prop:denseovert}
Let $(a_n)_{n \in \mathbb{N}}$ be a dense sequence in $\mathbf{X}$. Then $\nonempty_{\mathbf{X}}$ is computable relative to $(a_n)_{n \in \mathbb{N}}$.
\begin{proof}
If $U \in \mathcal{O}(\mathbf{X})$ is non-empty, it contains some $a_k$. Thus simultaneously testing $a_i \in U$ for all $i \in \mathbb{N}$ will detect non-emptyness of $U$, if true.
\end{proof}
\end{proposition}

As all represented spaces admit a dense sequence (by lifting a dense sequence in the domain of the representation), we see that all represented spaces are overt -- merely computable overtness survives as a distinguishing criterion.

\begin{proposition}
\label{prop:overthausdorff}
Let $\mathbf{X}$ be (computably) overt and $\mathbf{Y}$ be (computably) $T_2$. Then $\mathcal{C}(\mathbf{X}, \mathbf{Y})$ is (computably) $T_2$.
\begin{proof}
For $f, g : \mathbf{X} \to \mathbf{Y}$ we find $f \neq g \Leftrightarrow \nonempty_{\mathbf{X}}(\{x \in \mathbf{X} \mid f(x) \neq g(x)\})$.
\end{proof}
\end{proposition}

Similar to our procedure for compact sets, we can introduce the space $\mathcal{V}(\mathbf{X})$ of overt sets. The crucial mapping will be $\operatorname{Intersects} : \mathcal{V}(\mathbf{X}) \times \mathcal{O}(\mathbf{X}) \to \mathbb{S}$ defined by $\operatorname{Intersects}(A, U) = 1$ iff $A \cap U \neq \emptyset$. This makes clear that $\mathcal{V}(\mathbf{X})$ can be understood as a subspace of $\mathcal{O}(\mathcal{O}(\mathbf{X}))$. A set $O \in \mathcal{O}(\mathcal{O}(\mathbf{X}))$ will correspond to an element of $\mathcal{V}(\mathbf{X})$, iff it is of the form $O_A = \{U \in \mathcal{O}(\mathbf{X}) \mid U \cap A \neq \emptyset\}$ for some subset $A \subseteq X$. Note that $O_A = O_B$ for two sets $A, B \subseteq X$, iff $\overline{A} = \overline{B}$ (here $\overline{\phantom{A}}$ denotes the topological closure). Hence, we obtain $\mathcal{V}(\mathbf{X})$ by interpreting each set of the form $O_A \in \mathcal{O}(\mathcal{O}(\mathbf{X}))$ as $\overline{A} \subseteq X$.

In particular, we see that the overt sets and the closed sets coincide extensionally. In fact, the equivalence class of representations for $\mathcal{V}(\mathbf{X})$ has been studied for a long time, called the representation of closed sets by positive information. It is known that this representation is incomparable to the representation by negative information, i.e. that neither $\id : \mathcal{V}(\mathbf{X}) \to \mathcal{A}(\mathbf{X})$ nor $\id : \mathcal{A}(\mathbf{X}) \to \mathcal{V}(\mathbf{X})$ is continuous for non-empty $\mathbf{X}$. Our justification to avoid the word \emph{closed} for the space $\mathcal{V}(\mathbf{X})$ lies in the computable structure available on it, which differs significantly from the usual closure properties expected from closed sets:

\begin{proposition}
\label{prop:overtoperations}
Let $\mathbf{X}$, $\mathbf{Y}$, $\mathbf{Z}$ be represented spaces, and let $\mathbf{Z}$ be computably overt. Then the following functions are computable:
\begin{enumerate}
\item $\overline{\phantom{A}} : \mathcal{O}(\mathbf{Z}) \to \mathcal{V}(\mathbf{Z})$
\item $x \mapsto \overline{\{x\}} : \mathbf{X} \to \mathcal{V}(\mathbf{X})$
\item $\cup : \mathcal{V}(\mathbf{X}) \times \mathcal{V}(\mathbf{X}) \to \mathcal{V}(\mathbf{X})$
\item $\overline{\bigcup} : \mathcal{C}(\mathbb{N}, \mathcal{V}(\mathbf{X})) \to \mathcal{V}(\mathbf{X})$
\item $\overline{\cap} : \mathcal{V}(\mathbf{X}) \times \mathcal{O}(\mathbf{X}) \to \mathcal{V}(\mathbf{X})$
\item $\overline{\pi_1} : \mathcal{V}(\mathbf{X} \times \mathbf{Y}) \to \mathcal{V}(\mathbf{X})$
\item $(f, A) \mapsto \overline{f[A]} : \mathcal{C}(\mathbf{X}, \mathbf{Y}) \times \mathcal{V}(\mathbf{X}) \to \mathcal{V}(\mathbf{Y})$
\item $(f^{-1},A) \mapsto \overline{f[A]} :\subseteq \mathcal{C}(\mathcal{O}(\mathbf{Y}), \mathcal{O}(\mathbf{X})) \times \mathcal{V}(\mathbf{X}) \to \mathcal{V}(\mathbf{Y})$,  which takes a function between open sets that is the preimage map obtained from some (topologically continuous) function $f : \mathbf{X} \to \mathbf{Y}$ and an overt subset $A$ of $\mathbf{X}$, and outputs the closure of the image of $A$ under $f$ as an overt subset of $\mathbf{Y}$
\end{enumerate}
\begin{proof}
\begin{enumerate}
\item $\operatorname{Intersects}(\overline{V}, U) = \nonempty_{\mathbf{Z}}(U \cap V)$, using Proposition \ref{compana:prop:basicsetoperations} (3) for $\mathalpha{\cap} : \mathcal{O}(\mathbf{Z}) \times \mathcal{O}(\mathbf{Z}) \to \mathcal{O}(\mathbf{Z})$
\item $\operatorname{Intersects}(\overline{\{x\}}, U) = (x \in U)$, using Proposition \ref{compana:prop:basicsetoperations} (7) for $\mathalpha{\in} : \mathbf{X} \times \mathcal{O}(\mathbf{X}) \to \mathbb{S}$
\item $\operatorname{Intersects}(A \cup B, U) = \vee(\operatorname{Intersects}(A, U), \operatorname{Intersects}(B, U))$
\item $\operatorname{Intersects}(\overline{\bigcup_{n \in \mathbb{N}} A_n}, U) = \bigvee_{n \in \mathbb{N}}(\operatorname{Intersects}(A_n, U))$
\item $\operatorname{Intersects}(\overline{V \cap A}, U) = \operatorname{Intersects}(A, V \cap U)$, using Proposition \ref{compana:prop:basicsetoperations} (3) for $\mathalpha{\cap} : \mathcal{O}(\mathbf{X}) \times \mathcal{O}(\mathbf{X}) \to \mathcal{O}(\mathbf{X})$
\item $\operatorname{Intersects}(\overline{\pi_1(A)}, U) = \operatorname{Intersects}(A, U \times Y)$, using Proposition \ref{compana:prop:basicsetoperations} (8) for $\mathalpha{\times} : \mathcal{O}(\mathbf{X}) \times \mathcal{O}(\mathbf{Y}) \to \mathcal{O}(\mathbf{X} \times \mathbf{Y})$
\item $\operatorname{Intersects}(\overline{f[A]}, U) = \operatorname{Intersects}(A, f^{-1}(U))$, using Proposition \ref{compana:prop:basicsetoperations} (6) to obtain $f^{-1} : \mathcal{O}(\mathbf{Y}) \to \mathcal{O}(\mathbf{X})$
\item $\operatorname{Intersects}( \overline{f[A]}, U) = \operatorname{Intersects}(A, f^{-1}(U))$
\end{enumerate}
\end{proof}
\end{proposition}

\subsubsection*{Remarks}
The notion of overtness was named as such by \name{Taylor}, the definition here follows \name{Escard\'o} \cite{escardo}. This reference also contains Proposition \ref{prop:overthausdorff} and Proposition \ref{prop:overtoperations} (7). The statements of Proposition \ref{prop:denseovert} and \ref{prop:overtoperations} (5) are very similar to \cite[Proposition 2.2]{bauer4} by \name{Bauer} and \name{Lesnik} (2012). In restricted settings such as Euclidean spaces, \name{Weihrauch} has studied $\mathcal{V}(\mathbf{X})$ as the closed subsets of $\mathbf{X}$ represented by positive information, and e.g.~proven Proposition \ref{prop:overtoperations} (3) in \cite{weihrauchd}. The representation of $\mathcal{V}(\mathbf{X})$ via an identification with a subspace of $\mathcal{O}(\mathcal{O}(\mathbf{X}))$ was done by \name{Schr\"oder} in \cite[Section 4.4.2]{schroder5} (2002), together with several of its closure properties.

Both \name{Escard\'o} and \name{Taylor} have pointed out that classically, every space is overt.

\name{Collins} translated the definition of overtness into the language of computable analysis (under the name \emph{effectively separable}, and introduced the notation $\mathcal{V}(\mathbf{X})$ in \cite{collins5} (2010). He discusses the role of topological closure as the canonization operation for the overt sets; and provides the results of Proposition \ref{prop:overtoperations} ($1.-5.,7.$). His claimed characterization of the computably overt spaces as those having a computable dense sequence in \cite[Proposition 3.30]{collins5} (stated without proof) is wrong, though. A counterexample is the space $\{p \in \Cantor \mid p \textnormal{ is not computable}\}$ understood as a subspace of $\Cantor$. The failure of \cite[Proposition 3.30]{collins5} then impacts \cite[Theorem 3.32]{collins5}.

The space $\mathcal{V}(\mathbf{X})$ corresponds to the lower Vietoris topology (or, equivalently, the lower Fell topology) (as shown by \name{Schr\"oder}). This topology was introduced by \name{Vietoris} in 1922 \cite{vietoris}. Again, a general source for hyperspace topologies is \cite{beer}. Thus, we can obtain the Vietoris topology on a space of subsets as $\mathcal{K}(\mathbf{X}) \wedge \mathcal{V}(\mathbf{X})$ and the Fell topology as $\mathcal{A}(\mathbf{X}) \wedge \mathcal{V}(\mathbf{X})$. Both of these spaces appear in a variety of results.

\section{Discreteness}
Just as the $T_2$ separation property could be characterizes as making intersection of compact sets computable, there is a property of a space that makes intersection of overt sets computable. This property, \emph{discreteness}, is in many ways the dual to $T_2$ separation, just as overt is the dual to compact.

\begin{definition}
\label{def:discrete}
A represented space $\mathbf{X}$ is called (computably) discrete, iff the map $x \mapsto \{x\} : \mathbf{X} \to \mathcal{O}(\mathbf{X})$ is well-defined and continuous (computable).
\end{definition}

\begin{theorem}
\label{theo:discrete}
The following properties are equivalent for a represented space $\mathbf{X}$:
\begin{enumerate}
\item $\mathbf{X}$ is (computably) discrete.
\item $\id : \mathcal{V}(\mathbf{X}) \to \mathcal{O}(\mathbf{X})$ is well-defined and continuous (computable), and $\mathbf{X}$ is $T_1$.
\item $\cap : \mathcal{V}(\mathbf{X}) \times \mathcal{V}(\mathbf{X}) \to \mathcal{V}(\mathbf{X})$ is well-defined and continuous (computable), and $\mathbf{X}$ is $T_1$.
\item $\mathalpha{=} : \mathbf{X} \times \mathbf{X} \to \mathbb{S}$ defined by $\mathalpha{=}(x, x) = \top$ and $\mathalpha{=}(x, y) = \bot$ otherwise is continuous (computable).
\item $\Delta_\mathbf{X} = \{(x, x) \mid x \in \mathbf{X}\} \in \mathcal{O}(\mathbf{X} \times \mathbf{X})$ (is computable).
\end{enumerate}
\begin{proof}
\begin{description}
\item[$1. \Rightarrow 2.$] We find $x \in A$ iff $\operatorname{Intersects}(A, \{x\})$. The (topological) $T_1$ property is a straightforward consequence of singletons being open.
\item[$2. \Rightarrow 3.$] Use $\overline{\cap}$ from Proposition \ref{prop:overtoperations} (5) together with the assumption.
\item[$3. \Rightarrow 4.$] Given $(x, y) \in \mathbf{X} \times \mathbf{X}$, we can use Proposition \ref{prop:overtoperations} (2) to compute $(\overline{\{x\}}, \overline{\{y\}}) \in \mathcal{V}(\mathbf{X}) \times \mathcal{V}(\mathbf{X})$. The $T_1$ property means $\overline{\{x\}} = \{x\}$ and $\overline{\{y\}} = \{y\}$. Subsequently we use $\cap$ from the assumption to compute $\{x\} \cap \{y\} \in \mathcal{V}(\mathbf{X})$, and then $\operatorname{Intersects}(\{x\} \cap \{y\}, X)$, which is identical to $\mathalpha{=}(x,y)$.
\item[$4. \Rightarrow 1.$] This is a consequence of currying to $x \mapsto (y \mapsto \mathalpha{=}(x,y))$.
\item[$4. \Leftrightarrow 5.$] This is straightforward using $\mathcal{C}(\mathbf{X} \times \mathbf{X}, \mathbb{S}) \cong \mathcal{O}(\mathbf{X} \times \mathbf{X})$.
\end{description}
\end{proof}
\end{theorem}

As any represented space is separable (as pointed out in Section \ref{sec:overt}), the requirement that $\{x\} \in \mathcal{O}(\mathbf{X})$ for any $x \in \mathbf{X}$ already forces $\mathbf{X}$ to be countable, thus producing the following:

\begin{proposition}
\label{prop:discretecountable}
Any discrete represented space is countable.
\end{proposition}

It is worth pointed out that -- unlike in topology -- being computably discrete does not imply being computably Hausdorff. A counterexample can be constructed by starting with a recursively enumerable but not recursive set $A \subseteq \mathbb{N}$, and then identifying the elements of $A$. This means we define a representation $\delta : \Baire \to (\mathbb{N} \setminus A) \cup \{A\}$ by $\delta(p) = A$ if $p(0) \in A$ and $\delta(p) = p(0)$ if $p(0) \notin A$. This makes equality recognizable, but not refutable.

\begin{proposition}
\label{prop:compactdiscrete}
Let $\mathbf{X}$ be (computably) compact and $\mathbf{Y}$ be (computably) discrete. Then $\mathcal{C}(\mathbf{X}, \mathbf{Y})$ is (computably) discrete.
\begin{proof}
We find $\mathalpha{=}(f,g) = \operatorname{IsFull}_\mathbf{X}(\{x \in \mathbf{X} \mid \mathalpha{=}(f(x), g(x))\})$, the claim now follows with Theorem \ref{theo:discrete} (4) and Proposition \ref{compana:prop:compactnessbasics} (2).
\end{proof}
\end{proposition}

Unlike for its dual result in Proposition \ref{prop:overthausdorff}, it is easy to find a counter-example showing the necessity of a restriction on the domain in the preceding result. As such, consider $\mathbf{X} = \mathbf{Y} = \mathbb{N}$. The space $\mathbb{N}$ is computably discrete, but not compact, and $\mathcal{C}(\mathbb{N}, \mathbb{N}) \cong \Baire$ is not discrete by Proposition \ref{prop:discretecountable}.

\subsubsection*{Remarks}
\name{Escard\'o} uses the characterization in Theorem \ref{theo:discrete} (4) as definition of discreteness, and proceeds to prove Theorem \ref{theo:discrete} $[4. \Rightarrow 2.]$ and Proposition \ref{prop:compactdiscrete} in \cite{escardo}. That computable discreteness does not imply computably Hausdorff has been observed by both \name{Escard\'o}; and by \name{Weihrauch} in \cite{weihrauchh}.

In the setting of computable analysis, \name{Collins} introduced discreteness and proved the equivalence from Theorem \ref{theo:discrete} ($1. \Leftrightarrow 2.$) as \cite[Theorem 3.31 (1)]{collins5}.

\section{Admissibility as effective $T_0$ separation}
\label{section:admissibility}
Recall the computable maps $x \mapsto \sat{\{x\}} : \mathbf{X} \to \mathcal{K}(\mathbf{X})$ (Proposition \ref{prop:operationsoncompacts} (2)) and $x \mapsto \overline{\{x\}} : \mathbf{X} \to \mathcal{V}(\mathbf{X})$ (Proposition \ref{prop:overtoperations} (2)), and note that both are realized in the same way: An element $x \in \mathbf{X}$ provides the abstract capacity to determine membership in an open set (since $x \in U \Leftrightarrow \sat{\{x\}} \subseteq U \Leftrightarrow \overline{\{x\}} \cap U \neq \emptyset$).

\begin{definition}
Consider the computable map $\kappa_\mathbf{X} : \mathbf{X} \to \mathcal{O}(\mathcal{O}(\mathbf{X}))$ defined by $\kappa_\mathbf{X}(x) = \{U \in \mathcal{O}(\mathbf{X}) \mid x \in U\}$. Let $\mathbf{X}_\kappa$ denote the image of $\mathbf{X}$ in $\mathcal{O}(\mathcal{O}(\mathbf{X}))$ under $\kappa_\mathbf{X}$.
\end{definition}

With the considerations above, we see that we may consider $\mathbf{X}_\kappa$ simultaneously as a subspace of $\mathcal{K}(\mathbf{X})$ and $\mathcal{V}(\mathbf{X})$. If we understand $\mathbf{X}_\kappa$ as a subspace of $\mathcal{K}(\mathbf{X})$, then $\kappa_\mathbf{X}(x) = \sat{\{x\}}$, if as a subspace of $\mathcal{V}(\mathbf{X})$, then $\kappa_\mathbf{X}(x) = \overline{\{x\}}$. This identification will be helpful to see that $\kappa$ is an endofunctor\footnote{In fact, we could even call $\kappa$ a \emph{computable} endofunctor, as its restriction to any homset is computable.} on the category-extension of represented spaces. That $\kappa$ commutes with composition is rather obvious, hence we only need the following:

\begin{proposition}
\label{prop:kappafunctor}
There is an induced computable map $\kappa : \mathcal{C}(\mathbf{X}, \mathbf{Y}) \to \mathcal{C}(\mathbf{X}_\kappa,\mathbf{Y}_\kappa)$ for all represented spaces $\mathbf{X}$, $\mathbf{Y}$ such that the following diagram commutes:
$$\begin{CD}
\mathbf{X} @>f>> \mathbf{Y}\\
@VV\kappa_\mathbf{X}V @VV\kappa_\mathbf{Y}V\\
\mathbf{X}_\kappa @>\kappa(f)>> \mathbf{Y}_\kappa
\end{CD}$$
\begin{proof}
This follows e.g.~from Proposition \ref{prop:operationsoncompacts} (6) and $\mathbf{X}_\kappa \subseteq \mathcal{K}(\mathbf{X})$, together with $\sat{f[\sat{\{x\}}]} = \sat{\{f(x)\}}$.
\end{proof}
\end{proposition}

In general, $\kappa_\mathbf{X} : \mathbf{X} \to \mathbf{X}_\kappa$ will fail to be (computably) continuously invertible, however, those spaces admitting a continuous (computable) left-inverse for $\kappa_\mathbf{X}$ can be characterized as exactly those fully understandable in terms of their topology.

\begin{definition}
\label{def:admissible}
Call $\mathbf{X}$ (computably) admissible, iff $\kappa_\mathbf{X} : \mathbf{X} \to \mathbf{X}_\kappa$ admits a continuous (computable) left-inverse.
\end{definition}

We will proceed to demonstrate that $\kappa$ is the coreflector of the (cartesian-closed) subcategory of the (computably) admissible represented spaces inside the category of represented spaces. In a way, this coreflector is analogous to the Kolmogorov-quotient (or $T_0$-coreflector) in topology.
\begin{corollary}
\label{corr:admissiblereflection}
Let $\mathbf{X}$, $\mathbf{Y}$ be represented spaces, and let $\mathbf{Y}$ be (computably) admissible. There is a (computable) continuous map $\mathfrak{R} : \mathcal{C}(\mathbf{X}, \mathbf{Y}) \to \mathcal{C}(\mathbf{X}_\kappa, \mathbf{Y})$ such that $f = \mathfrak{R}(f) \circ \kappa_\mathbf{X}$ for all $f \in \mathcal{C}(\mathbf{X}, \mathbf{Y})$.
\begin{proof}
Take the map $\kappa : \mathcal{C}(\mathbf{X}, \mathbf{Y}) \to \mathcal{C}(\mathbf{X}_\kappa, \mathbf{Y}_\kappa)$ from Proposition \ref{prop:kappafunctor}. Then we use Proposition \ref{compana:prop:functionspacesbasics} (4, 6) to compose $\kappa(f)$ with $\kappa_\mathbf{Y}^{-1}$ obtained from Definition \ref{def:admissible}.
\end{proof}
\end{corollary}

\begin{proposition}
\label{prop:sierpadmissible}
$\mathbb{S}$ is computably admissible.
\begin{proof}
We can explicitly define $\kappa_\mathbb{S}^{-1} :\subseteq \mathcal{O}(\mathcal{O}(\mathbb{S})) \to \mathbb{S}$ by $\kappa_\mathbb{S}^{-1}(U) = (\{\top\} \in U)$.
\end{proof}
\end{proposition}

\begin{theorem}
\label{theo:admissiblecc}
Let $\mathbf{Y}$ be (computably) admissible. Then $\mathcal{C}(\mathbf{X}, \mathbf{Y})$ is (computably) admissible.
\begin{proof}
Given some $U \in \mathcal{O}(\mathbf{Y})$ and $x \in \mathbf{X}$, we can compute $\{f \in \mathcal{C}(\mathbf{X}, \mathbf{Y}) \mid f(x) \in U\} \in \mathcal{O}(\mathcal{C}(\mathbf{X}, \mathbf{Y}))$. Hence, given $x \in \mathbf{X}$ and $\{U \in \mathcal{O}(\mathcal{C}(\mathbf{X}, \mathbf{Y})) \mid f \in U\} \in \mathcal{O}(\mathcal{O}(\mathcal{C}(\mathbf{X}, \mathbf{Y})))$ we can compute $\{U \in \mathcal{O}(\mathbf{Y}) \mid f(x) \in U\} \in \mathcal{O}(\mathcal{O}(\mathbf{Y}))$ by Proposition \ref{compana:prop:basicsetoperations} (6). By assumption, the latter suffices to compute $f(x) \in \mathbf{Y}$. Currying yields the claim.
\end{proof}
\end{theorem}

\begin{corollary}
$\mathbf{X}_\kappa$, $\mathcal{O}(\mathbf{X})$, $\mathcal{A}(\mathbf{X})$, $\mathcal{K}(\mathbf{X})$ and $\mathcal{V}(\mathbf{X})$ are all computably admissible.
\end{corollary}

\begin{corollary}
\label{corr:topologicalnotionskappa}
$\mathcal{O}(\mathbf{X}) \cong \mathcal{O}(\mathbf{X}_\kappa)$, $\mathcal{A}(\mathbf{X}) \cong \mathcal{A}(\mathbf{X}_\kappa)$, $\mathcal{K}(\mathbf{X}) \cong \mathcal{K}(\mathbf{X}_\kappa)$ and $\mathcal{V}(\mathbf{X}) \cong \mathcal{V}(\mathbf{X}_\kappa)$.
\end{corollary}

\begin{corollary}
\label{corr:kappaopen}
$\kappa_\mathbf{X}$ is effectively open, i.e.~there is a well-defined and computable map $K : \mathcal{O}(\mathbf{X}) \to \mathcal{O}(\mathbf{X}_\kappa)$ such that $K(U) = \kappa_{\mathbf{X}}[U]$.
\end{corollary}

\begin{corollary}
$\mathbf{X}$ is (computably) compact, (computably) $T_2$, (computably) overt or (computably) discrete if and only if $\mathbf{X}_\kappa$ has that property.
\end{corollary}

\begin{theorem}
\label{theo:admissiblecharac}
The following properties are equivalent for a represented space $\mathbf{X}$:
\begin{enumerate}
\item $\mathbf{X}$ is (computably) admissible.
\item $f \mapsto f_{\textrm{K}} : \mathcal{C}(\mathbf{Y}, \mathbf{X}) \to \mathcal{C}(\mathcal{K}(\mathbf{Y}), \mathcal{K}(\mathbf{X}))$ has a well-defined and continuous (computable) partial inverse for any represented space $\mathbf{Y}$; where $f_{\textrm{K}}(A) = \uparrow f[A]$.
\item $f \mapsto f^{-1} : \mathcal{C}(\mathbf{Y}, \mathbf{X}) \to \mathcal{C}(\mathcal{O}(\mathbf{X}), \mathcal{O}(\mathbf{Y}))$ has a well-defined and continuous (computable) partial inverse for any represented space $\mathbf{Y}$.
\item Any topologically continuous function $f : \mathbf{Y} \to \mathbf{X}$ (i.e.~$f^{-1} : \mathcal{O}(\mathbf{X}) \to \mathcal{O}(\mathbf{Y})$ is well-defined) is continuous as a function between represented spaces (i.e.~$f \in \mathcal{C}(\mathbf{Y}, \mathbf{X})$). (no computable counterpart)
\end{enumerate}
\begin{proof}
\begin{description}
\item[$1. \Rightarrow 2.$] For precision, let $f_K$ be the input, i.e.~the map from compact sets to compact sets derived from the desired output $f$. Now note $f = \kappa_\mathbf{X}^{-1} \circ f_K \circ \kappa_\mathbf{Y}$, and consider Proposition \ref{prop:kappafunctor}, Definition \ref{def:admissible} and Proposition \ref{compana:prop:functionspacesbasics} (4).
\item[$2. \Rightarrow 3.$] Given a continuous function of the form $f^{-1} : \mathcal{O}(\mathbf{X}) \to \mathcal{O}(\mathbf{Y})$ induced by some continuous $f : \mathbf{Y} \to \mathbf{X}$, we may use Proposition \ref{prop:operationsoncompacts} (4) to obtain the induced function $f_\textrm{K} : \mathcal{K}(\mathbf{Y}) \to \mathcal{K}(\mathbf{X})$. By assumption, we can recover $f \in \mathcal{C}(\mathbf{X}, \mathbf{Y})$ from the latter.
\item[$3. \Rightarrow 1.$]   First we shall show that $\kappa_\mathbf{X}^{-1} : \mathbf{X}_\kappa \to \mathbf{X}$ is well-defined, using proof-by-contradiction. So let us assume $\kappa_\mathbf{X}$ were not injective, i.e.~ there were $x \neq y \in \mathbf{X}$ with $\kappa_\mathbf{X}(x) = \kappa_\mathbf{X}(y)$, hence $x \in U \Leftrightarrow y \in U$ for any $U \in \mathcal{O}(\mathbf{X})$. Then the constant function $\overline{x}, \overline{y} : \{0\} \to \mathbf{X}$ with $\overline{x}(0) = x$ and $\overline{y}(0) = y$ are continuous and distinct, but still $\overline{x}^{-1} = \overline{y}^{-1}$ as maps of the type $\mathcal{O}(\mathbf{X}) \to \mathcal{O}(\{0\})$. This contradicts the well-definedness of $f^{-1} \mapsto f :\subseteq \mathcal{C}(\mathcal{O}(\mathbf{X}), \mathcal{O}(\{0\})) \to \mathcal{C}(\{0\}, \mathbf{X})$.

    Consider $\mathbf{Y} := \mathbf{X}_\kappa$. By Corollary \ref{corr:kappaopen}, $\kappa_{\mathbf{X}}$ is effectively open, hence $\kappa_\mathbf{X} \in \mathcal{C}(\mathcal{O}(\mathbf{X}), \mathcal{O}(\mathbf{X}_\kappa))$ is a computable element. As $\kappa_\mathbf{X}^{-1}$ is well-defined, $\kappa_\mathbf{X} = (\kappa_\mathbf{X}^{-1})^{-1}$, and we can use the assumption to obtain $\kappa_\mathbf{X}^{-1} : \mathbf{X}_\kappa \to \mathbf{X}$ from $\kappa_\mathbf{X}$ as open map.
\item[$1. \Rightarrow 4.$] By Corollary \ref{corr:topologicalnotionskappa}, a function $f : \mathbf{Y} \to \mathbf{X}$ is topologically continuous, iff it is so as a function $f : \mathbf{Y} \to \mathbf{X}_\kappa$. We may curry topologically continuous functions, too, and obtain $\chi_f : \mathbf{Y} \times \mathcal{O}(\mathbf{X}) \to \mathbb{S}$ as a topologically continuous function defined by $\chi_f(y, U) = \top$ iff $f(y) \in U$. By composition, so is $\chi_f \circ (\delta_\mathbf{Y} \times \delta_{\mathcal{O}(\mathbf{X})}) : \subseteq \Cantor \to \mathbb{S}$. But this just means $\chi_f \circ (\delta_\mathbf{Y} \times \delta_{\mathcal{O}(\mathbf{X})}) \in \mathcal{O}(\dom(\delta_\mathbf{Y} \times \delta_{\mathcal{O}(\mathbf{X})}))$. Now the process can be reversed, using represented space continuity in place of topological continuity to find that $f \in \mathcal{C}(\mathbf{Y}, \mathbf{X}_\kappa)$. Admissibility of $\mathbf{X}$, i.e.~$\mathbf{X} \cong \mathbf{X}_\kappa$, then allows us the make the final step and reach $f \in \mathcal{C}(\mathbf{Y}, \mathbf{X})$.
\item[$4. \Rightarrow 1.$] Again, we first show that $\kappa_\mathbf{X}^{-1}$ is well-defined using proof-by-contradiction. Assume there were $x \neq y \in \mathbf{X}$ with $\kappa_\mathbf{X}(x) = \kappa_\mathbf{X}(y)$, hence $x \in U \Leftrightarrow y \in U$ for any $U \in \mathcal{O}(\mathbf{X})$. Then any function $f : \Cantor \to \{x, y\} \subseteq \mathbf{X}$ is topologically continuous. However, there are $2^{2^{\aleph_0}}$ such functions, whereas $|\mathcal{C}(\Cantor, \mathbf{X})| \leq 2^{\aleph_0}$.

    Using again $\kappa_\mathbf{X} = (\kappa_\mathbf{X}^{-1})^{-1}$, as well as Corollary \ref{corr:kappaopen}, we see that $\kappa_\mathbf{X}^{-1} : \mathbf{X}_\kappa \to \mathbf{X}$ is topologically continuous, hence we obtain it as a continuous function between represented spaces.
\end{description}
\end{proof}
\end{theorem}

Note that only the equivalence of $1., 2., 3.$ in the preceding theorem has a proper place in a synthetic treatment -- the reference in $4.$ to the mere (external) well-definedness of $f^{-1}$ is meaningless from a strictly internal view on a given category. Its provability here is due to the definition of continuity for represented spaces via (topological) continuity on Cantor space. What we obtain from it is the observation that the admissibly represented spaces (with continuous maps) simultaneously form a subcategory of the represented spaces and of the topological spaces, moreover, that this is in some sense the largest such joint subcategory.

As an example for the interplay of admissibility and some of the other properties of represented spaces, we shall briefly revisit the connection between a function and its graph (cf.~Proposition \ref{compana:prop:t1characterization} (7)):

\begin{proposition}
\label{prop:graphinverse1}
Let $\mathbf{Y}$ be (computably) admissible and (computably) compact. Then $\operatorname{Graph}^{-1} : \subseteq \mathcal{A}(\mathbf{X} \times \mathbf{Y}) \to \mathcal{C}(\mathbf{X}, \mathbf{Y})$ is continuous (computable), where $\dom(\operatorname{Graph}^{-1}) = \{A \in \mathcal{A}(\mathbf{X} \times \mathbf{Y}) \mid \exists f : \mathbf{X} \to \mathbf{Y} \  \ \operatorname{Graph}(f) = A\}$.
\begin{proof}
We may assume that $\operatorname{Graph}(f) \in \mathcal{A}(\mathbf{X} \times \mathbf{Y})$ and $x \in \mathbf{X}$ are given and show how to obtain $f(x) \in \mathbf{Y}$. Using Proposition \ref{compana:prop:basicsetoperations} (1,9), we obtain $\{y \mid (x,y) \in \operatorname{Graph}(f)\} = \{f(x)\} \in \mathcal{A}(\mathbf{Y})$. Corollary \ref{compana:corr:idclosedcompact} gets us $\sat{\{f(x)\}} \in \mathcal{K}(\mathbf{Y})$, by definition of $\mathbf{Y}_\kappa$ thus $f(x) \in \mathbf{Y}_\kappa$. By definition of admissibility, this in turn suffices to obtain $f(x) \in \mathbf{Y}$.
\end{proof}
\end{proposition}

Similar results can be obtained for compact and overt graphs:
\begin{proposition}
\label{prop:graphinverse2}
Let $\mathbf{Y}$ be (computably) admissible and $\mathbf{X}$ (computably) $T_2$. Then $\operatorname{Graph}^{-1} : \subseteq \mathcal{K}(\mathbf{X} \times \mathbf{Y}) \to \mathcal{C}(\mathbf{X}, \mathbf{Y})$ is continuous (computable), where $\dom(\operatorname{Graph}^{-1}) = \{A \in \mathcal{K}(\mathbf{X} \times \mathbf{Y}) \mid \exists f : \mathbf{X} \to \mathbf{Y} \  \ \operatorname{Graph}(f) = A\}$.
\begin{proof}
We may assume that $\operatorname{Graph}(f) \in \mathcal{K}(\mathbf{X} \times \mathbf{Y})$ and $x \in \mathbf{X}$ are given and show how to obtain $f(x) \in \mathbf{Y}$. By Definition \ref{compana:def:t1} we can obtain $\{x\} \in \mathcal{A}(\mathbf{X})$, and then $\{x\} \times Y \in \mathcal{A}(\mathbf{X} \times \mathbf{Y})$ by Proposition \ref{compana:prop:basicsetoperations} (8). Then Proposition \ref{prop:operationsoncompacts} (4) provides us with $\sat{\left ((\{x\} \times Y) \cap \operatorname{Graph}(f)\right )} = \sat{(\{x\} \times \{f(x)\})} \in \mathcal{K}(\mathbf{X} \times \mathbf{Y})$. Then we use projection (Proposition \ref{prop:operationsoncompacts} (9)) to get $\sat{\{f(x)\}} \in \mathcal{K}(\mathbf{Y})$, admissibility of $\mathbf{Y}$ enables us to extract $f(x)$.
\end{proof}
\end{proposition}

\begin{proposition}
\label{prop:graphinverse3}
Let $\mathbf{Y}$ be (computably) admissible and $\mathbf{X}$ (computably) discrete. Then $\operatorname{Graph}^{-1} : \subseteq \mathcal{V}(\mathbf{X} \times \mathbf{Y}) \to \mathcal{C}(\mathbf{X}, \mathbf{Y})$ is continuous (computable), where $\dom(\operatorname{Graph}^{-1}) = \{A \in \mathcal{V}(\mathbf{X} \times \mathbf{Y}) \mid \exists f : \mathbf{X} \to \mathbf{Y} \  \ \operatorname{Graph}(f) = A\}$.
\begin{proof}
We may assume that $\operatorname{Graph}(f) \in \mathcal{V}(\mathbf{X} \times \mathbf{Y})$ and $x \in \mathbf{X}$ are given and show how to obtain $f(x) \in \mathbf{Y}$. By Definition \ref{def:discrete} we can obtain $\{x\} \in \mathcal{O}(\mathbf{X})$, and then $\{x\} \times Y \in \mathcal{O}(\mathbf{X} \times \mathbf{Y})$ by Proposition \ref{compana:prop:basicsetoperations} (8). Then Proposition \ref{prop:overtoperations} (5) provides us with $\overline{(\{x\} \times Y) \cap \operatorname{Graph}(f)} = \overline{\{x\} \times \{f(x)\}} \in \mathcal{V}(\mathbf{X} \times \mathbf{Y})$. Then we use projection (Proposition \ref{prop:overtoperations} (6)) to get $\overline{\{f(x)\}} \in \mathcal{V}(\mathbf{Y})$, admissibility of $\mathbf{Y}$ enables us to extract $f(x)$.
\end{proof}
\end{proposition}

\subsubsection*{How to use admissibility}
Admissibility is at the core of a very common scheme to prove computability of some specific mapping:

\begin{enumerate}
\item Conclude that the solution set $S$ is closed based on its definition (this often involves the Hausdorff condition).
\item Obtain some compact candidate set $K$ (either from the situation, or by assumption - this often involves some bounds).
\item Compute $\sat{(S \cap K)}$ as a compact set via Proposition \ref{prop:operationsoncompacts} (4).
\item Use domain-specific reasoning or assumption to conclude that the solution $s$ is unique.
\item As we have $\sat{\{s\}}$ available as a compact set, if our target space is admissible, we can compute $s$.
\end{enumerate}

While this scheme is usually well-hidden and obscured, it is present e.g.~in \cite{hoyrup} by \name{Galatolo}, \name{Hoyrup} and \name{Rojas}, \cite{rettinger} by \name{Rettinger}, \cite{collins6,collins6b} by \name{Collins} and \name{Gra\c{c}a}. A very similar algorithm (albeit in a slightly different formal setting) is described by \name{Escard\'o} in \cite{escardo6}. This use of admissibility also underlies the model of \emph{non-deterministic type-2 machines} suggest by \name{Ziegler} \cite{ziegler7} and studied further by \name{Brattka}, \name{de Brecht} and the author in \cite{paulybrattka,paulybrattka2}.

\subsubsection*{Remarks}
The notion of admissibility was introduced by \name{Kreitz} and \name{Weihrauch} in 1985 \cite{kreitz}. They essentially define a represented space to admissible if it is isomorphic to a subspace of $\mathcal{O}(\mathbb{N})$, which captures the countably based admissible spaces. Then they proceed to prove that this implies the condition in Theorem \ref{theo:admissiblecharac} (4).

The understanding of admissibility as presented here is essentially due to \name{Schr\"oder}. In \cite{schroder} (2002) \name{Schr\"oder} uses the condition in Theorem \ref{theo:admissiblecharac} (4) as definition of admissibility, characterizes the topologies arising as $\mathcal{O}(\mathbf{X})$ and in particular proves a slightly weaker version of Theorem \ref{theo:admissiblecc} (as it requires admissibility of $\mathbf{X}$, too). In his thesis \cite{schroder5} (2002), \name{Schr\"oder} also defines computable admissibility (i.e.~Definition \ref{def:admissible}) and provides the statements of Propositions \ref{prop:kappafunctor}, \ref{prop:sierpadmissible}, Corollary \ref{corr:admissiblereflection} and Theorems \ref{theo:admissiblecc}, \ref{theo:admissiblecharac}. \name{Lietz} also added to the development, e.g.~\cite[Theorem 3.2.7.]{lietz}, which shows that Definition \ref{def:admissible} is indeed well-suited for its purpose.

The identification of a point $x \in \mathbf{X}$ with its neighbourhood filter $\{U \in \mathcal{O}(\mathbf{X}) \mid x \in U\} \in \mathcal{O}(\mathcal{O}(\mathbf{X}))$ is reminiscent of the ultrafilter approach to topology.

The study of admissibility can be seen as an attempt to generalize the coincide of computability relative to an oracle and (topological) continuity beyond the setting of Baire space. While dealing only with a restricted class of spaces, \cite{brattkahertling2} by \name{Brattka} and \name{Hertling} (1994) and \cite{paulyziegler} by the author and \name{Ziegler} (2013) consider the questions for multivalued functions rather than just functions.

The connections between a continuous function and its graph as exemplified in Propositions \ref{prop:graphinverse1}, \ref{prop:graphinverse2}, \ref{prop:graphinverse3} has been advanced by \name{Brattka} in \cite{brattka6} (2008).

\section{Open relations and compact or overt sets}
The dual nature of the spaces of overt and of compact sets becomes clearer when the interaction with open relations is studied. This in particular generalizes the well-known results about upper and lower computability of the maximum of compact subsets of the real numbers presented either as compact sets, or as overt sets (i.e.~by positive information).

\begin{proposition}
\label{prop:exists}
The map $\exists : \mathcal{O}(\mathbf{X} \times \mathbf{Y}) \times \mathcal{V}(\mathbf{X}) \to \mathcal{O}(\mathbf{Y})$ defined by $\exists(R, A) = \{y \in Y \mid \exists x \in A \ (x, y) \in R\}$ is computable. Moreover, whenever $\exists : \mathcal{O}(\mathbf{X} \times \mathbf{Y}) \times \mathcal{S}(\mathbf{X}) \to \mathcal{O}(\mathbf{Y})$ is computable for some hyperspace $\mathcal{S}(\mathbf{X})$ and some space $\mathbf{Y}$ containing a computable element $y_0$, then $\overline{\phantom{A}} : \mathcal{S}(\mathbf{X}) \to \mathcal{V}(\mathbf{X})$ is computable.
\begin{proof}
To see that $\exists$ is computable, note that $y \in \exists(R, A)$ iff $\operatorname{Intersects}(A, \operatorname{Cut}(y, R))$ and use Proposition \ref{compana:prop:basicsetoperations} (9) and the definition of $\mathcal{V}(\mathbf{X})$. Now if $\exists$ is computable for some other hyperspace $\mathcal{S}(\mathbf{X})$ in place of $\mathcal{V}(\mathbf{X})$, we can compute $\overline{\phantom{A}}: \mathcal{S}(\mathbf{X}) \to \mathcal{V}(\mathbf{X})$ via $\operatorname{Intersects}(\overline{A}, U) = \mathalpha{\in}(y_0,\exists(U \times Y, A))$.
\end{proof}
\end{proposition}

\begin{corollary}
\label{corr:unionovert}
$\bigcup : \mathcal{V}(\mathcal{O}(\mathbf{Y})) \to \mathcal{O}(\mathbf{Y})$ is computable. (This functions maps $A \in \mathcal{V}(\mathcal{O}(\mathbf{Y}))$ to $\left (\bigcup_{U \in A} U \right ) \in \mathcal{O}(\mathbf{Y})$.)
\begin{proof}
Pick $\mathbf{X} = \mathcal{O}(\mathbf{Y})$, and instantiate $R$ with $\{(U, x) \in \mathcal{O}(\mathbf{Y}) \times \mathbf{Y} \mid x \in U\}$.
\end{proof}
\end{corollary}

\begin{proposition}
\label{prop:forall}
The map $\forall : \mathcal{O}(\mathbf{X} \times \mathbf{Y}) \times \mathcal{K}(\mathbf{X}) \to \mathcal{O}(\mathbf{Y})$ defined by $\forall(R, A) = \{y \in Y \mid \forall x \in A \ (x, y) \in R\}$ is computable. Moreover, whenever $\forall : \mathcal{O}(\mathbf{X} \times \mathbf{Y}) \times \mathcal{S}(\mathbf{X}) \to \mathcal{O}(\mathbf{Y})$ is computable for some hyperspace $\mathcal{S}(\mathbf{X})$ and some space $\mathbf{Y}$ containing a computable element $y_0$, then $\sat{\id} : \mathcal{S}(\mathbf{X}) \to \mathcal{K}(\mathbf{X})$ is computable.
\begin{proof}
To see that $\forall$ is computable, note that $y \in \forall(R, A)$ iff $\operatorname{IsContainedIn}(A, \operatorname{Cut}(y, R))$ and use Proposition \ref{compana:prop:basicsetoperations} (9) and Proposition \ref{prop:operationsoncompacts} (1). Now if $\forall$ is computable for some other hyperspace $\mathcal{S}(\mathbf{X})$ in place of $\mathcal{K}(\mathbf{X})$, we can compute $\sat{\id}: \mathcal{S}(\mathbf{X}) \to \mathcal{K}(\mathbf{X})$ via $\operatorname{IsContainedIn}(A, U) = \mathalpha{\in}(y_0,\forall(U \times Y, A))$.
\end{proof}
\end{proposition}

\begin{corollary}
\label{corr:intersectioncompact}
$\bigcap : \mathcal{K}(\mathcal{O}(\mathbf{Y})) \to \mathcal{O}(\mathbf{Y})$ is computable. (This function maps $A \in \mathcal{K}(\mathcal{O}(\mathbf{Y}))$ to $\left (\bigcap_{U \in A} U \right ) \in \mathcal{O}(\mathbf{Y})$.)
\begin{proof}
Pick $\mathbf{X} = \mathcal{O}(\mathbf{Y})$, and instantiate $R$ with $\{(U, x) \in \mathcal{O}(\mathbf{Y}) \times \mathbf{Y} \mid x \in U\}$.
\end{proof}
\end{corollary}

To make the connection to the computability of maxima on the reals, we first generalize the Dedekind-cut construction of the reals. Let $\mathalpha{\prec} \in \mathcal{O}(\mathbf{X} \times \mathbf{X})$ be some transitive and open relation. Now we define $\overline{\mathbf{X}_\prec} := \{U \in \mathcal{O}(\mathbf{X}) \mid \forall x \in \mathbf{X} \left (\exists y \in U \ x \prec y \Rightarrow x \in U\right )\}$ as the subspace of $\mathcal{O}(\mathbf{X})$ containing the open initial segments of $\prec$. Next, consider the computable map $x \mapsto \operatorname{Cut}(x, \mathalpha{\prec}) : \mathbf{X} \to \overline{\mathbf{X}_\prec}$ mapping $x$ to $\{y \in \mathbf{X} \mid y \prec x\}$. We shall denote its image by $\mathbf{X}_\prec$, and identify an element $x$ with its $\prec$-lower ideal. Depending on the properties of $\prec$, the map $\id : \mathbf{X} \to \mathbf{X}_\prec$ induced by such an identification may fail to be injective. Spaces of this form that have been studied so far are $\mathbb{R}_<$ and $\mathbb{R}_>$, the former in particular playing a central r\^ole in computable measure theory (e.g.~\name{Schr\"oder} \cite{schroder2}, \name{Collins} \cite{collins4}). By construction, the spaces $\mathbf{X}_\prec$ and $\overline{\mathbf{X}_\prec}$ are computably admissible.

\begin{proposition}
\label{prop:dedekindbasics}
\begin{enumerate}
\item $\prec \in \mathcal{O}(\mathbf{X} \times \mathbf{X}_\prec)$, $\prec \in \mathcal{O}(\mathbf{X}_\succ \times \mathbf{X})$ are computable.
\item $(x, y) \mapsto \{z \in \mathbf{X} \mid x \prec z \prec y\} : \mathbf{X}_\succ \times \mathbf{X}_\prec \to \mathcal{O}(\mathbf{X})$ is computable.
\item If $\prec$ is dense\footnote{Here, \emph{dense} is used in the order-theoretic sense, i.e.~$\prec$ is dense means that $\forall x, y \in \mathbf{X} x \prec y \Rightarrow \left ( \exists z \in \mathbf{X} \ x \prec z \prec y \right )$.} and $\mathbf{X}$ is (computably) overt, then $\prec \in \mathcal{O}(\mathbf{X}_\succ \times \mathbf{X}_\prec)$ (is computable).
\item $\sup_\prec : \mathcal{V}(\mathbf{X}) \to \overline{\mathbf{X}_\prec}$ and $\sup_\prec : \mathcal{K}(\mathbf{X}) \to \overline{\mathbf{X}_\succ}$ are computable.
\item Let $\mathbf{Y} \subseteq \mathbf{X}$ be a dense\footnote{And here \emph{dense} is used in the topological sense, i.e.~means$\forall U \in \mathcal{O}(\mathbf{X}) U \neq \emptyset \rightarrow \exists x \in \mathbf{Y} \cap U$.} and computably overt subspace, and $\prec$ be dense. Then we can identify $\overline{\mathbf{Y}_\prec}$ and $\overline{\mathbf{X}_\prec}$.
\end{enumerate}
\begin{proof}
\begin{enumerate}
\item $\prec : \mathbf{X} \times \mathbf{X}_\prec \to \mathbb{S}$ is a restriction of the computable map $\in : \mathbf{X} \times \mathcal{O}(\mathbf{X}) \to \mathbb{S}$ from Proposition \ref{compana:prop:basicsetoperations} (7). The second claim follows by symmetry.
\item This map is a restriction of the computable map $\cap : \mathcal{O}(\mathbf{X}) \times \mathcal{O}(\mathbf{X}) \to \mathcal{O}\mathbf{X})$ from Proposition \ref{compana:prop:basicsetoperations} (3).
\item If $\prec$ is dense, then $x \prec y$ iff $\nonempty_\mathbf{X}(\{z \in \mathbf{X} \mid x \prec z \prec y\})$ which is computable by 2. and Definition \ref{def:overtspace}.
\item These are corollaries of Propositions \ref{prop:exists}, \ref{prop:forall}.
\item The restriction of any element of $\overline{\mathbf{X}_\prec}$ to $\mathbf{Y}$ will yield an element of $\overline{\mathbf{Y}_\prec}$. Under the conditions given, we can recover a lower ideal $U \in \overline{\mathbf{X}_\prec}$ from its restriction to $\mathbf{Y}$ by noting $x \in U \Leftrightarrow \exists y \in \mathbf{Y} \ x \prec y \wedge y \in U$.
\end{enumerate}
\end{proof}
\end{proposition}

If $\prec$ is e.g.~a linear order, then the identification $\id : \mathbf{X} \to \mathbf{X}_\prec$ is injective, however, will usually not be computably invertible. The best we can obtain here is the following:

\begin{proposition}
\label{prop:intervalbase}
Let the intervals of $\prec$ be an effective base for $\mathcal{O}(\mathbf{X})$, i.e.~let the computable map $(a_i, b_i)_{i \in \mathbb{N}} \mapsto \bigcup_{i \in \mathbb{N}} \{x \mid a_i \prec x \prec b_i\} : \mathcal{C}(\mathbb{N}, \mathbf{X} \times \mathbf{X}) \to \mathcal{O}(\mathbf{X})$ be surjective and admit a computable right-inverse (as a multivalued function). Then $\mathbf{X}_\kappa \cong (\mathbf{X}_\prec \wedge \mathbf{X}_\succ)$.
\begin{proof}
As $\mathbf{X}_\prec$ and $\mathbf{X}_\succ$ are admissible, the identification maps $\id : \mathbf{X}_\kappa \to \mathbf{X}_\prec$ and $\id : \mathbf{X}_\kappa \to \mathbf{X}_\succ$ are well-defined and computable. It only remains to be shown that they admit a joint computable right-inverse. For this, we need to demonstrate how given $x \in \mathbf{X}_\prec \wedge \mathbf{X}_\succ$ and $U \in \mathcal{O}(\mathbf{X})$ we can recognize $x \in U?$. By assumption, $U$ can be effectively expressed as $\bigcup_{i \in \mathbb{N}} \{y \mid a_i \prec y \prec b_i\}$, and by Proposition \ref{prop:dedekindbasics} (1) we can simultaneously semidecide $a_i \prec x$ and $x \prec b_i$, until for some $i \in \mathbb{N}$ both statements are true, in which case $x \in U$ is accepted.
\end{proof}
\end{proposition}

We just briefly introduce the represented space $\mathbb{R}$ to provide an actual example. In this, let $\nu_\mathbb{Q}$ be a standard notation of the rationals. Then we define $\rho_\mathbb{R} : \subseteq \Baire \to \mathbb{R}$ by $\rho_\mathbb{R}(p) = x$ iff $\forall n \in \mathbb{N} \ . \ d(\nu_\mathbb{Q}(p(n)), x) < 2^{-n}$; and consider $\rho_\mathbb{R}$ as the representation of $\mathbb{R}$. We find $< \in \mathcal{O}(\mathbb{R} \times \mathbb{R})$ to be computable, thus $\mathbb{R}$ is computably $T_2$. Moreover, $\mathbb{R}$ satisfies the conditions of Proposition \ref{prop:intervalbase}.

\begin{corollary}
$\overline{\mathbb{R}_\prec} \cong \overline{\mathbb{Q}_\prec}$, $\overline{\mathbb{R}_\succ} \cong \overline{\mathbb{Q}_\succ}$, $\mathbb{R} \cong \mathbb{R}_\prec \wedge \mathbb{R}_\succ$.
\end{corollary}

\begin{corollary}
\label{corr:maxcomp}
$\max : \mathcal{K}(\mathbb{R}) \wedge \mathcal{V}(\mathbb{R}) \to \mathbb{R}$ is computable.
\end{corollary}

\begin{corollary}
$\operatorname{max-value} : (\mathcal{K}(\mathbb{R}) \wedge \mathcal{V}(\mathbb{R})) \times \mathcal{C}(\mathbb{R}, \mathbb{R}) \to \mathbb{R}$ is computable, where $\operatorname{max-value}(A, f) = \max \{f(x) \mid x \in A\}$.
\begin{proof}
By Proposition \ref{prop:operationsoncompacts} (7) and Proposition \ref{prop:overtoperations} (7), we can compute $f[A] \in (\mathcal{K}(\mathbb{R}) \wedge \mathcal{V}(\mathbb{R}))$, then we use $\max$ from Corollary \ref{corr:maxcomp}.
\end{proof}
\end{corollary}

\subsubsection*{Remarks}
Corollaries \ref{corr:unionovert} and \ref{corr:intersectioncompact} are present in \name{Escard\'o}'s \cite{escardo} (2004). The topological counterpart of Corollary \ref{corr:intersectioncompact} has been described by \name{Nachbin} \cite{nachbin} in 1992. The results on $\mathbb{R}_<$ and $\mathbb{R}_>$ are e.g.~included in \name{Weihrauch}'s \cite{weihrauchd} (2000).

\section{Concluding remarks}
It was demonstrated that compactness (Proposition \ref{compana:prop:compactnessbasics}), $T_2$ separation (Proposition \ref{compana:prop:t1characterization}), overtness and discreteness (to a lesser extent), and admissibility (Theorem \ref{theo:admissiblecharac}) all encompass various equivalent properties of represented spaces each characterized by the continuity or computability of specific maps. This marks the extent of possible generalizations of many standard results on computability for sets and functions. In particular, unnecessary restrictions in prior work such as admissibility, second-countability or metrizability are avoided.

As the realizers witnessing the computability of the relevant mappings are generic and independent of the represented spaces involved, the approach presented here is a viable alternative to the introduction of multi-representations of countably based admissibly represented spaces suggested in \cite{rettinger2} by \name{Rettinger} and \name{Weihrauch}.

Moving away from effective topological spaces, i.e.~countably based admissible represented spaces, to represented spaces as the primitive objects of investigation for set and function computability opens up options for an integration with the study of hyper-computation. \name{Ziegler} \cite{ziegler2, ziegler3} and \name{Brattka} (see \cite{paulybrattka}) observed that various kinds of hyper-computation such as limit computability or computability with mindchanges can be adequately characterized by means of operators generating new represented spaces from given ones. Based on these, a suitable characterization of functions such as topological closure and interior in terms of computable maps between appropriate represented spaces ought to be possible along the lines of \cite{gherardi3}.

A more general approach to such \emph{jump operators} can be found in \cite{debrecht5} by \name{de Brecht}. These operators (which amount to endofunctors) form the basis of a new approach to descriptive set theory \cite{paulydebrecht, pauly-descriptive} as suggested by \name{Pauly} and \name{de Brecht}, which in particular provides derived represented spaces e.g.~for the $\Sigma_n^0$-measurable functions between given represented spaces, thus generalizing \cite{brattka} by \name{Brattka}. In this, the applicability of descriptive set theory is extended even further than to the Quasi-Polish spaces introduced by \name{de Brecht} in \cite{debrecht6}.

The abstract study of represented spaces seems to hold further aspects to contemplate. Proposition \ref{prop:intervalbase} already made use of the concept of an effective basis, and a systematic study of those seems to be promising also in order to obtain results on product spaces, e.g.~a version of the Tychonoff theorem. Computable measure theory also is amenable to a similar development as computable topology underwent here, as demonstrated by \name{Collins} in \cite{collins4}.
\bibliographystyle{eptcs}
\bibliography{../spieltheorie}

\section*{Acknowledgements}
I would like to thank Vasco Brattka, Martin Escard\'o and Matthias Schr\"oder for various discussions crucial for the formation of the present work.

The article benefitted from comments by Eike Neumann, Greg Yang and Klaus Weihrauch, as well as the anonymous referees.
\end{document}